\documentclass[12pt,reqno]{amsart}

\usepackage{amssymb}
\usepackage[arrow]{xypic}

\DeclareMathOperator{\chr}{char}\DeclareMathOperator{\decc}{dec}
\DeclareMathOperator{\Fix}{Fix} \DeclareMathOperator{\Gal}{Gal}
 \DeclareMathOperator{\im}{im}
\DeclareMathOperator{\Ker}{ker} \DeclareMathOperator{\proj}{proj}
\DeclareMathOperator{\res}{res} \DeclareMathOperator{\tra}{tra}

\newcommand{\dec}{{\decc}}
\newcommand{\Ec}{\mathcal{E}}
\newcommand{\F}{\mathbb{F}}
\newcommand{\Fp}{\F_p}
\newcommand{\Ft}{\F_2}
\newcommand{\G}{{\Gamma}}
\newcommand{\GG}{\mathcal{G}}
\renewcommand{\H}{\mathcal{H}}

\newcommand{\Ic}{\mathcal{I}}
\newcommand{\Jc}{\mathcal{J}}
\newcommand{\Kc}{\mathcal{K}}
\newcommand{\N}{{\Delta}}
\newcommand{\Nb}{\mathbb{N}}
\newcommand{\Q}{\mathbb{Q}}
\newcommand{\Z}{\mathbb{Z}}

\begin{document}

\title[Detecting Pro-$p$-groups]{Detecting Pro-$p$-groups that
are not Absolute Galois Groups, Expanded Version}

\author[Dave Benson]{Dave Benson}
\address{Department of Mathematical Sciences, University of
Aberdeen, Meston Building, King's College, Aberdeen AB24 3UE,
Scotland UK}
\email{bensondj@maths.abdn.ac.uk}

\author[Nicole Lemire]{Nicole Lemire}
\author[J\'{a}n Min\'{a}\v{c}]{J\'an Min\'a\v{c}}
\address{Department of Mathematics, Middlesex College, \ University
of Western Ontario, London, Ontario \ N6A 5B7 CANADA}
\email{nlemire@uwo.ca}
\email{minac@uwo.ca}

\author[John Swallow]{John Swallow}
\address{Department of Mathematics, Davidson College, Box 7046,
Davidson, North Carolina \ 28035-7046 USA}
\email{joswallow@davidson.edu}

\date{October 20, 2006}

\thanks{Nicole Lemire was supported in part by Natural Sciences
and Engineering Research Council of Canada grant R3276A01. J\'an
Min\'a\v{c} was supported in part by Natural Sciences and
Engineering Research Council of Canada grant R0370A01.  John Swallow
was supported in part by National Science Foundation grant
DMS-0600122.}

\maketitle

\def\thepart{A}

\newtheorem{theorem}{Theorem}[part]
\newtheorem{proposition}{Proposition}[part]
\newtheorem{corollary}[theorem]{Corollary}
\newtheorem*{corollary*}{Corollary}
\newtheorem{lemma}{Lemma}[part]

\theoremstyle{definition}
\newtheorem*{remark*}{Remark}
\newtheorem*{example*}{Example}

\parskip 10pt plus 2pt

Let $p$ be a prime.  It is a fundamental problem to classify the
absolute Galois groups $G_F$ of fields $F$ containing a primitive
$p$th root of unity.  In this paper we present several constraints
on such $G_F$, using restrictions on the cohomology of index $p$
normal subgroups from \cite{LMS}. In section~\ref{se:t} we classify
all maximal $p$-elementary abelian-by-order $p$ quotients of these
$G_F$.  In the case $p>2$, each such quotient contains a unique
closed index $p$ elementary abelian subgroup.  This seems to be the
first case in which one can completely classify nontrivial quotients
of absolute Galois groups by characteristic subgroups of normal
subgroups.  In section~\ref{se:asb} we derive analogues of theorems
of Artin-Schreier and Becker for order $p$ elements of certain small
quotients of $G_F$. Finally, in sections~\ref{se:fam}--\ref{se:fam3}
we construct new families of pro-$p$-groups which are not absolute
Galois groups over any field $F$.

As a consequence of our results, we prove the following limitations
on relator shapes of pro-$p$ absolute Galois groups. For elements
$\sigma$ and $\tau$ of a group $\G$, let ${}^0[\sigma,\tau]=\tau$,
${}^1[\sigma,\tau]=\sigma\tau\sigma^{-1}\tau^{-1}$, and
${}^n[\sigma,\tau] = [\sigma, {}^{n-1}[\sigma,\tau]]$ for $n\ge 2$.
Similarly, for subsets $\G_1$ and $\G_2$ of $\G$, let
${}^n[\G_1,\G_2]$ denote the closed subgroup generated by all
elements of the form ${}^n[\gamma_1, \gamma_2]$ for $\gamma_i\in
\G_i$.
\begin{theorem}\label{th:1}
    Let $p$ be an odd prime, $\G$ a pro-$p$-group with maximal
    closed subgroup $\N$, and $\sigma\in \G\setminus \N$.
    \begin{enumerate}
        \item Suppose that for some $\tau\in \N$ and some $2\le e\le
        p-2$
        \begin{align*}
            {}^e[\sigma,\tau]&\not\in {}^{p-1}[\sigma,\N]\Phi(\N) \\
            {}^{e+1}[\sigma,\tau] &\in \Phi(\N).
        \end{align*}
        Then $\G$ is not an absolute Galois group.

        \

        \noindent Moreover, if $\G$ contains a normal subgroup
        $\Lambda\subset \N$ such that $\G/\Lambda\simeq \Z/p^2\Z$,
        we may take $1\le e\le p-2$.
        \item Suppose that for some $\tau_1, \tau_2\in \N$,
        \begin{align*}
            [\sigma,\tau_i]&\not\in {}^{p-1}[\sigma,\N]\Phi(\N),
            &i=1, 2 \\ {}^2[\sigma,\tau_i] &\in \Phi(\N), &i=1,
            2\\ \langle [\sigma,\tau_1] \rangle \Phi(\N) &\neq
            \langle [\sigma,\tau_2] \rangle \Phi(\N)
        \end{align*}
        Then $\G$ is not an absolute Galois group.
        \item Suppose that
        \begin{equation*}
            \sigma^p \in {}^{2}[\sigma,\N]\Phi(\N).
        \end{equation*}
        Then $\G$ is not an absolute Galois group.
    \end{enumerate}
\end{theorem}
\noindent Here $\Phi(\N)=\N^p[\N,\N]$ denotes the Frattini subgroup
of $\N$.

A striking corollary is that one-relator pro-$p$-groups with
relations quite similar to the relations of Demu\v{s}kin groups for
odd primes \cite{D} cannot be absolute Galois groups.
\begin{corollary*}\label{co:1}
    Let $p$ be an odd prime and $\G$ a pro-$p$-group minimally
    generated by $\{\sigma_1,\sigma_2\} \cup \{\sigma_i\}_{i\in
    \Ic}$ subject to a single relation
    \begin{equation*}
        \sigma_1^q\cdot {}^f[\sigma_1,\sigma_2]\cdot
        \prod_{(i,j)\in \Jc} [\sigma_i,\sigma_j]\cdot
        \prod_{k \in \Kc} [\sigma_1^p,\sigma_k]
    \end{equation*}
    for some $2\le f\le p-1$, $q\in \Nb\cup \{0\}$ with $q = 0 \bmod
    (p^2)$, a finite ordered set of pairs $\Jc\subset \Ic\times
    \Ic$, and a finite ordered subset $\Kc$ of $\Ic$.  Then $\G$ is
    not an absolute Galois group.
\end{corollary*}

The results in \cite{LMS} may be used to establish further new
results on possible $V$-groups of fields and metabelian quotients of
absolute Galois pro-$p$-groups. (For the definition of the $V$-group
$V_F$ of a field $F$, see section~\ref{se:asb}.) Moreover, some of
the results here hold in a greater generality than their
formulations here.  For instance, the examples here of pro-$p$
non-absolute Galois groups are also examples of groups which are not
maximal pro-$p$-quotients of absolute Galois groups, by \cite[\S
6]{LMS}.  Furthermore, pro-$p$-groups which are not absolute Galois
groups are not $p$-Sylow subgroups of absolute Galois groups.  We
plan a systematic study of these concerns in \cite{BeMS}.

We observe that because this paper is concerned only with degree $1$
and degree $2$ cohomology, the results cited from \cite{LMS} rely
only on the Merkurjev-Suslin Theorem \cite[Theorem 11.5]{MeSu} and
not the recent work on the full Bloch-Kato Conjecture.

\section{$T$-groups}\label{se:t}

A \emph{$T$-group} is a nontrivial pro-$p$-group $T$ with a maximal
closed subgroup $N$ that is abelian of exponent dividing $p$. Then
$N$ is a normal subgroup, and the factor group $T/N$ acts naturally
on $N$: choose a lift $t\in T$ and act via $n\mapsto tnt^{-1}$.
Given any profinite group $\G$ with a closed normal subgroup $\N$ of
index $p$, the factor group $\G/\N^p[\N,\N]$ is a $T$-group.

Now suppose that $E/F$ is a cyclic field extension of degree $p$.
We define the \emph{$T$-group of $E/F$} to be $T_{E/F} :=
G_F/G_E^p[G_E,G_E]$.  In this section we classify those $T$-groups
realizable as $T_{E/F}$ for fields $F$ containing a primitive $p$th
root of unity or of characteristic $p$.

We develop a complete set of invariants $t_i$, $i=1, 2, \dots, p$,
and $u$ of $T$-groups as follows. For a pro-$p$-group $\G$, we
denote by $\G^{(n)}$ the $n$th group in the $p$-central series of
$\G$ and $\G_{(n)}$ the $n$th group in the central series of $\G$.
Hence $\G^{(1)}=\G_{(1)}=\G$, and, for $n\in \Nb$,
$\G^{(n+1)}=(\G^{(n)})^p\cdot [\G,\G^{(n)}]$ and
$\G_{(n+1)}=[\G,\G_{(n)}]$.  For a $T$-group $T$ we define
\begin{align*}
    t_1 &= \dim_{\Fp} H^1\left(\frac{Z(T)[p]}{Z(T)\cap T_{(2)}},
    \Fp\right) \\ t_i &= \dim_{\Fp} H^1 \left(\frac{{Z(T)\cap
    T_{(i)}}}{Z(T)\cap T_{(i+1)}},\Fp\right), \quad 2\le i\le p\\
    u &= \max\{i: 1\le i\le p, \ T^p\subset T_{(i)}\}.
\end{align*}
For readability we write $t_i$ and $u$ instead of the more
cumbersome $t_i(T)$ and $u(T)$.  Here $Z(T)$ denotes the center of
$T$ and $Z(T)[p]$ the subgroup of $Z(T)$ of elements of order
dividing $p$.  From our definition of $T$-groups it follows that
$T_{(2)} = [T,N] \subset N$ and therefore $Z(T) \cap T_{(2)} \subset
Z(T)[p]$, so that $t_1$ is well-defined.

\begin{proposition}\label{pr:tc}
    For arbitrary cardinalities $t_i$, $i=1, 2, \dots, p$, and
    $u$ with $1\le u\le p$, the following are equivalent:
    \begin{enumerate}
        \item\label{it:tc1} The $t_i$ and $u$ are invariants of some
        $T$-group
        \item\label{it:tc2}
        \begin{enumerate}
            \item\label{it:tc2a} if $u<p$ then $t_u\ge 1$, and
            \item\label{it:tc2b} if $u=p$ and $t_i=0$ for all $2\le
            i\le p$, then $t_1\ge 1$.
        \end{enumerate}
    \end{enumerate}
    Moreover, $T$-groups are uniquely determined up to isomorphism
    by these invariants.
\end{proposition}

\begin{remark*}
    In condition (\ref{it:tc2b}) above, one can omit the condition
    $u=p$: from (\ref{it:tc2a}) it follows that if $u<p$ and $t_i =
    0$ for all $2 \le i \le p$, then $u=1$ and therefore $t_1 \ge
    1$.
\end{remark*}

In the following theorem, we consider fields $F$ either containing a
primitive $p$th root of unity or of characteristic $p$.
\begin{theorem}\label{th:t} Suppose $p$ is an odd prime. The
    following are equivalent:
    \begin{enumerate}
        \item $T$ is a $T$-group with invariants $t_i$ and $u$
        satisfying
        \begin{enumerate}
            \item $u\in \{1,2\}$,
            \item $t_2=u-1$, and
            \item $t_i=0$ for $3\le i<p$.
        \end{enumerate}
        \item $T\simeq T_{E/F}$ for some cyclic extension of degree
        $p$ of a field $F$ as above.
    \end{enumerate}

    Suppose $p=2$.  Then the following are equivalent:
    \begin{enumerate}
        \item $T$ is a $T$-group
        \item $T\simeq T_{E/F}$ for some cyclic extension of degree
        $2$ of a field $F$.
    \end{enumerate}
\end{theorem}

\begin{corollary*}
    Suppose $p$ is an odd prime, $T$ is a finite nonabelian
    $T$-group, $N$ is a maximal closed subgroup of $T$ that is
    abelian of exponent dividing $p$, and $u$ is the invariant
    of $T$ defined above.  Set $e=\dim_{\Fp} H^1(N,\Fp)$ and
    $f=\dim_{\Fp} Z(T)$. If
    \begin{equation*}
        e \not\equiv f+(u-1) \mod p-1,
    \end{equation*}
    then $T\not\simeq T_{E/F}$ for any cyclic extension $E/F$
    of degree $p$.
\end{corollary*}

Let $G$ be a cyclic group of order $p$ and $M_i$, $i=1, 2, \dots,
p$, denote the unique cyclic $\Fp G$-module of dimension $i$.  Since
the $\Fp G$-modules we consider will be multiplicative groups, we
usually write the action of $G$ exponentially. Recall that for any
generator $g$ of $G$, the element $(g-1)\in \Fp G$ is nilpotent of
degree $p$ and the fixed $\Fp G$-submodule $M_i^G =
M_i^{(g-1)^{i-1}}$ for $1\le i\le p$. For a set $\Ic$, let
$M_i^{\Ic}$ denote the product of $\vert\Ic\vert$ copies of $M_i$
endowed with the product topology.

We use the word \emph{duality} exclusively to refer to Pontrjagin
duality, between compact and discrete abelian groups and more
generally between compact and discrete $\Fp G$-modules.  Since we
shall consider only the case when $G$ has order $p$, we do not need
to pass to completed group rings.  (See, for instance, \cite[\S
7.1--7.3 and \S 14]{Koc} and \cite{Br}.) We denote the Pontrjagin
dual of a discrete or compact group $\Gamma$ by $\Gamma^\vee$.

\begin{lemma}\label{le:t}
    Suppose $T$ is a $T$-group with invariants $t_i$, $i=1, 2,
    \dots, p$ and $u$, and $N$ is a maximal closed subgroup of $T$
    that is abelian of exponent dividing $p$.  Let $\sigma\in
    T\setminus N$ and set $\bar\sigma$ to be the image of $\sigma$
    in $G:=T/N$.  Then
    \begin{enumerate}
        \item\label{it:lt1} $N$ is a topological $\Fp G$-module.
        \item\label{it:lt2} For any $\Fp G$-submodule $M$ of $N$ and
        $i\ge 0$, ${}^i[T,M]=M^{(\bar\sigma-1)^i}$.
        \item\label{it:lt3} For $i\ge 2$, $T_{(i)} =
        N^{(\bar\sigma-1)^{i-1}}$.
        \item\label{it:lt4} There exist sets $\Ic_i$, $i=1, 2,
        \dots, p$, such that $N$ decomposes into indecomposable $\Fp
        G$-modules as $N = M_1^{\Ic_1}\times M_2^{\Ic_2} \times
        \cdots \times M_p^{\Ic_p}$, endowed with the product
        topology, and all such decompositions of $N$ are equivalent.
        Moreover, for $i\ge 2$, $t_i=\vert \Ic_i\vert$, and
        \begin{equation*}
            t_1 = \begin{cases}
                1+\vert \Ic_1\vert, &T \text{\ is abelian\ of\
                exponent\ } p\\
                \vert \Ic_1\vert, &\text{otherwise}.
            \end{cases}
        \end{equation*}
        \item\label{it:lt5} $T^p = \langle \sigma^p \rangle \cdot
        T_{(p)}$.
        \item\label{it:lt6} If $u<p$ then $u$ is the minimal $v$ with
        $1\le v\le p-1$ such that there is a commutative diagram of
        pro-$p$-groups
        \begin{equation*}
            \xymatrix{ 1 \ar[r]& N \ar[r] \ar@{>>}[d] & T \ar[r]
            \ar@{>>}[d] & G \ar[r] \ar[d]^{=} & 1 \\ 1 \ar[r] & M_v
            \ar[r] & H \ar[r] & G \ar[r] & 1}
        \end{equation*}
        with a lift of $\bar\sigma$ in $H$ of order $p^2$. If $u=p$
        then no such diagram exists for $1\le v\le p-1$.
    \end{enumerate}
\end{lemma}

\begin{proof}
    Since $N$ is a maximal closed subgroup of a pro-$p$-group $T$,
    $N$ is normal of index $p$.  Therefore $G$ is well-defined and
    is a cyclic group of order $p$ with generator $\bar\sigma$.

    (\ref{it:lt1}).  Since $N$ is abelian, the action of $T$ on $N$
    factors through $G$. This action on the pro-$p$-group $N$ is
    continuous, whence $N$ is a topological $\Fp G$-module.

    (\ref{it:lt2}).  Let $M$ be an arbitrary $\Fp G$-submodule of
    $N$. Suppose that $\tau\in T$ is arbitrary, and write $\tau=
    n\sigma^{i}$ for $n\in N$ and $i \in \Nb \cup \{0\}$.  Since the
    action of $T$ factors through $G$ we see that $[\tau,M] =
    M^{(\bar\sigma^i-1)}$.  Since $(\bar\sigma-1)$ divides
    $(\bar\sigma^i-1)$ for all $i\in \Nb$ we obtain that $[T,M]=
    M^{(\bar\sigma-1)}$, and the result follows by induction.
    (Moreover, we observe that if $p$ does not divide $i$, then
    $[\tau,M] =[M,\tau]=[T,M]=[M,T]$.)

    (\ref{it:lt3}).  Since $N$ is abelian of index $p$, see that
    $[T,T]=[T,N]$, and the result follows from \eqref{it:lt2}.

    (\ref{it:lt4}).  Because $\Fp G$ is an Artinian principal ideal
    ring, every $\Fp G$-module $U$ decomposes into a direct sum of
    cyclic $\Fp G$-modules \cite[Theorem~6.7]{SV}.  Every cyclic
    $\Fp G$-module is indecomposable, self-dual, and local.
    Moreover, each such module has a local endomorphism ring.  By
    the Krull-Schmidt-Azumaya Theorem (see \cite[Theorem 12.6]{AF}),
    all decompositions of an $\Fp G$-module $U$ into indecomposables
    are equivalent. (In our special case one can check this fact
    directly.)  Applying these results to $U=N^{\vee}$ and using
    duality (see \cite[Lemma 2.9.4 and Theorem 2.9.6]{RZ}), we obtain
    the decomposition and its uniqueness.

    Observe that for an $\Fp G$-submodule $M$ of $N$, $Z(T)\cap M =
    M^G$.  Using \eqref{it:lt3} together with $M_j^G =
    M_j^{(\bar\sigma-1)^{j-1}}$ for all $1\le j\le p$, we calculate
    \begin{equation*}
        Z(T)\cap N = \left(M_1\right)^{\Ic_i} \times
        \left(M_2^{(\bar\sigma-1)}\right)^{\Ic_2}\times \cdots
        \times \left(M_p^{(\bar\sigma-1)^{p-1}}\right)^{\Ic_p}.
    \end{equation*}
    and for $i\ge 2$
    \begin{equation*}
        Z(T)\cap T_{(i)} = \begin{cases}
        \left(M_i^{(\bar\sigma-1)^{i-1}}\right)^{\Ic_i} \times
        \cdots \times
        \left(M_p^{(\bar\sigma-1)^{p-1}}\right)^{\Ic_p}, &2\le i\le
        p \\ \{1\}, &i>p. \end{cases}
    \end{equation*}
    Then from $\dim_{\Fp} M_j^{(\bar\sigma-1)^{j-1}}=1$, $1\le j\le
    p$, we deduce that $t_i=\vert \Ic_i\vert$, $2\le i\le p$.

    For the case $i=1$, suppose first that $T$ is abelian of
    exponent $p$.  Then $Z(T)=Z(T)[p]=T$ and $T_{(2)}=\{1\}$.  By
    \eqref{it:lt3}, $N^{(\bar\sigma-1)}=\{1\}$, from which we deduce
    $\vert \Ic_i \vert=0$ for $i\ge 2$.  Therefore
    \begin{align*}
        t_1 &= \dim_{\Fp} H^1(T,\Fp) = \dim_{\Fp} H^1(T/N,\Fp) +
        \dim_{\Fp} H^1(N,\Fp) \\ &= 1 + \vert \Ic_1 \vert,
    \end{align*}
    as desired.

    Now suppose that $T$ is nonabelian.  If an element $\delta
    \in T \setminus N$ commutes with $N$, then $T$ is abelian.
    Hence no element of $T\setminus N$ commutes with $N$ and thus
    $Z(T)\subset N$.  We obtain $Z(T)[p]=Z(T)\cap N$ and so
    \begin{equation*}
        t_1 = \dim_{\Fp} H^1 \left(\frac{Z(T)\cap N}{Z(T)\cap
        T_{(2)}}, \Fp \right) = \vert\Ic_1\vert.
    \end{equation*}

    Finally, assume that $T$ is abelian and not of exponent $p$.
    Then the endomorphism $\phi:x\mapsto x^p$ satisfies $\im \phi
    \neq \{1\}$ and $N\subset \ker \phi$. Since $N$ is a maximal
    closed subgroup, we obtain $\ker \phi = N$, that is, $N=T[p]$.
    Since $Z(T)=T$, we deduce $t_1 = \dim_{\Fp} H^1(N,\Fp) = \vert
    \Ic_1 \vert$.

    (\ref{it:lt5}). For $\delta\in N$ we have $(\delta\sigma)^2 =
    [\sigma,\delta]\delta^2\sigma^2$ and, by induction,
    \begin{equation*}
        (\delta\sigma)^i =
        {[\underbrace{\sigma,[\sigma,\dots,[\sigma}_{i-1\ \text{times}},
        \delta]\cdots]]}^{\Tiny\begin{pmatrix}i\\ i\end{pmatrix}} \cdots
        [\sigma,[\sigma,\delta]]^{\Tiny\begin{pmatrix}
        i\\ 3\end{pmatrix}} [\sigma,\delta]^{\Tiny\begin{pmatrix}i\\
        2\end{pmatrix}} \delta^i\sigma^i.
    \end{equation*}
    Then $(N\sigma)^p= (\sigma^p)\cdot {}^{p-1}[\sigma, N]$, which
    by \eqref{it:lt2} and \eqref{it:lt3} may be written $\sigma^p
    \cdot T_{(p)}$. Replacing $\sigma$ with $\sigma^v$ for
    $(v,p)=1$, we conclude $T^p = \langle \sigma^p \rangle \cdot
    T_{(p)}$.

    (\ref{it:lt6}). Suppose that for some $v<u$ there is a commutative
    diagram
    \begin{equation*}
        \xymatrix{ 1 \ar[r]& N \ar[r] \ar@{>>}[d] & T \ar[r]
        \ar@{>>}[d] & G \ar[r] \ar[d]^{=} & 1 \\ 1 \ar[r] & M_v
        \ar[r] & H \ar[r] & G \ar[r] & 1}
    \end{equation*}
    with a lift of $\bar\sigma$ in $H$ of order $p^2$.  Then the
    surjection $T\twoheadrightarrow H$ factors through
    $T/N^{(\bar\sigma-1)^{v}}$. But by \eqref{it:lt3}, $T_{(v+1)} =
    N^{(\bar\sigma-1)^v}$ and by definition of $u$, we have
    $T^p\subset T_{(u)}\subset T_{(v+1)}$. Hence every lift of
    $\bar\sigma$ into $T/N^{(\bar\sigma-1)^v}$ is of order $p$, and
    the same holds for $H$.  We conclude that no commutative diagram
    as above with $\bar\sigma$ lifting to an element of order $p^2$
    exists for $v<u$.

    Going the other direction, suppose that $u<p$ and consider
    $\sigma^p$.  By \eqref{it:lt5}, $T^p=\langle \sigma^p\rangle
    \cdot T_{(p)}$, and $T_{(p)}\subset T_{(u+1)}$.  By definition
    of $u$, we have $T^p\subset T_{(u)}$ and $T^p\not\subset
    T_{(u+1)}$. We conclude that $\sigma^p\in T_{(u)}\setminus
    T_{(u+1)}$. Now $\sigma^p\in N$ since $N$ is a maximal closed
    subgroup. By \eqref{it:lt3}, since
    $T_{(u)}=N^{(\sigma-1)^{u-1}}$ for $u\ge 2$, we obtain
    $\sigma^p\in N^{(\sigma-1)^{u-1}}\setminus N^{(\sigma-1)^u}$.
    From $[\sigma,\sigma^p]=1$ we obtain $\sigma^p\in N^G$.
    Therefore
    \begin{equation*}
        \sigma^p \in N^G\cap N^{(\sigma-1)^{u-1}} \setminus
        N^G \cap N^{(\sigma-1)^u}.
    \end{equation*}

    We claim that there exists an $\Fp G$-submodule $M_u$ of $N$
    such that $M_u^G=\langle \sigma^p\rangle$ and $N=M_u\times
    \tilde N$ for some $\Fp G$-submodule $\tilde N$ of $N$. Assume a
    factorization of $N$ into cyclic $\Fp G$-submodules as in
    \eqref{it:lt4}.  Let $w \in N$ such that $w^{(\sigma-1)^{u-1}} =
    \sigma^p$.  Consider the components of $w$ in the factors of
    $N$.  Without loss of generality the nonzero components of $w$
    lie in factors of dimension at least $u$, and there must exist
    at least one such factor $M_{u}$ of dimension precisely $u$ such
    that the projection $\proj_{M_{u}} w$ of $w$ generates $M_{u}$
    as an $\Fp G$-module.  Let $\tilde N$ be the $\Fp G$-submodule
    of $N$ generated by all of the factors of $N$ except $M_{u}$.
    Then it follows that $N = M_u \times \tilde N$.  Factoring $T$
    by $\tilde N$, we obtain the commutative diagram
    \begin{equation*}
        \xymatrix{ 1 \ar[r]& N \ar[r] \ar@{>>}[d] & T \ar[r]
        \ar@{>>}[d] & G \ar[r] \ar[d]^{=} & 1 \\ 1 \ar[r] &
        M_u\ar[r] & T/\tilde N \ar[r] & G \ar[r] & 1}
    \end{equation*}
    in which a lift of $\bar\sigma$ is of order $p^2$.

    We have shown that if $u<p$, then $u$ is the minimal $v$ in
    $1\le v\le p-1$ such that a commutative diagram
    \begin{equation*}
        \xymatrix{ 1 \ar[r]& N \ar[r] \ar@{>>}[d] & T \ar[r]
        \ar@{>>}[d] & G \ar[r] \ar[d]^{=} & 1 \\ 1 \ar[r] &
        M_v\ar[r] & H \ar[r] & G \ar[r] & 1}
    \end{equation*}
    exists with a lift of $\bar\sigma$ of order $p^2$, and if $u=p$
    then no such commutative diagram with a lift of $\bar\sigma$ of
    order $p^2$ exists for $1\le v\le p-1$.
\end{proof}

\begin{remark*}
    Observe that from Lemma~\ref{le:t}(\ref{it:tb3},\ref{it:tb4})
    it follows that the nilpotent index of $T$ is at most $p$.
    Moreover, this index is $\le p-1$ if and only if $\Ic_p =
    \varnothing$.  In this case, each finite quotient of $T$ is a
    regular $p$-group.  (See \cite[\S 12.4]{Ha}.) This fact can also
    be seen from the identity expressing $(\sigma\delta)^p$ as a
    product of $p$th-powers and commutators as in the proof of
    \eqref{it:lt5} above.
\end{remark*}

\begin{lemma}\label{le:t1}
    Let $T$ be a $T$-group with invariants $t_i$, $i=1, 2, \dots, p$,
    and $u$.
    \begin{enumerate}
        \item\label{it:tb1} $T$ is abelian if and only if $t_i=0$
        for all $2\le i\le p$.
        \item\label{it:tb2} $T$ is of exponent $p$ if and only if
        $u=p$ and $t_u=0$.
        \item\label{it:tb3} If $u<p$ then $t_u\ge 1$.
        \item\label{it:tb4} If $t_i=0$ for all $2\le i\le
        p$, then $t_1\ge 1$.
    \end{enumerate}
\end{lemma}

\begin{proof}
    (\ref{it:tb1}).  If $T$ is abelian then $T_{(2)}=\{1\}$ and
    therefore $t_i=0$ for $i\ge 2$.  Conversely, suppose that
    $t_i=0$ for $i\ge 2$.  Then by Lemma~\ref{le:t}\eqref{it:lt4},
    we have $N^G = N$.  By Lemma~\ref{le:t}\eqref{it:lt3},
    $T_{(2)}=\{1\}$, whence $T$ is abelian.

    (\ref{it:tb2}).  If $u=p$ and $t_p=0$ then the fact that $T$ is
    of exponent $p$ follows from the definitions.  We show the
    converse, as follows.  Assume that $T$ is of exponent $p$. From
    Lemma~\ref{le:t}(\ref{it:lt5}) we have $T_{(p)}= T^p=\{1\}$.
    Hence $t_p=0$, and the definition of $u$ gives $u=p$.

    (\ref{it:tb3}).  Suppose that $u<p$.  Since $T^p\subset T_{(u)}$
    and $T^p\not\subset T_{(u+1)}$, we cannot have $T^p=\{1\}$.  Let
    $\delta\in T$ be arbitrary with $\delta^p\neq 1$. Then
    $\delta\not\in N$, and since $N$ is of index $p$, $\delta^p\in
    N$.  Since $[\delta,\delta^p]=1$ we deduce $\delta^p\in Z(T)\cap
    N$.  We have obtained $T^p\subset Z(T)\cap T_{(u)}$ and
    $T^p\not\subset Z(T)\cap T_{(u+1)}$, whence $t_u\ge 1$.

    (\ref{it:tb4}).  Suppose $t_i=0$ for all $2\le i\le p$.  From
    \eqref{it:tb1} we conclude that $T$ is abelian.  Hence from the
    definition of $t_1$ and the assumption that $T$ is nontrivial,
    we conclude that $t_1 \ge 1$.
\end{proof}

\begin{remark*}
    Observe that if $p=2$ then (\ref{it:tb2}) implies (\ref{it:tb1})
    in Lemma~\ref{le:t1}, which is good because every group of
    exponent $2$ is abelian!
\end{remark*}

\begin{proof}[Proof of Proposition~\ref{pr:tc}]
    By Lemma~\ref{le:t1}(\ref{it:tb3},\ref{it:tb4}), we see that
    conditions \eqref{it:tc2} on the $t_i$ and $u$ hold for any
    $T$-group.  Hence \eqref{it:tc1} implies \eqref{it:tc2}.  We now
    show that \eqref{it:tc2} implies \eqref{it:tc1}.  Suppose we are
    given cardinalities $t_i$, $i=1, 2, \dots, p$, and $u$
    satisfying conditions~\eqref{it:tc2}.

    \emph{The case $u<p$}. Let $G$ be a group of order $p$ and
    $\Ic_i$, $i=1, \cdots, p$ sets with cardinalities $\vert
    \Ic_i\vert$ satisfying $\vert \Ic_i\vert = t_i$ for $i\neq u$
    and $1+\vert \Ic_u\vert = t_u$. Set $N = X \times M_1^{\Ic_1}
    \times M_2^{\Ic_2} \times \cdots \times M_p^{\Ic_p}$, where
    $X\simeq M_u$ and $N$ is endowed with the product topology.
    Observe that $N$ possesses a system of neighborhoods of the
    identity consisting of open $G$-invariant normal subgroups.
    Hence $N$ is a pseudocompact $\Fp G$-module. (See
    \cite[page~443]{Br}.)

    Define an action of $\Z_p$ on $N$ by letting a generator
    $\sigma$ of $\Z_p$ act via a generator of $G$. (We write $\Z_p$
    multiplicatively to reduce confusion.)  We see that $\Z_p$ acts
    continuously on $N$ and we form the semidirect product $N\rtimes
    \Z_p$ in the category of pro-$p$-groups. Now choose an $\Fp
    G$-module generator $x$ of $X$ and define $x_i=x^{(\sigma-1)^i}$
    for $0\le i\le u$.  Since $(\sigma-1)$ is nilpotent of degree
    $u$ on $X$ we obtain $x_{u-1}\neq 1$ and $x_u=1$.  We set $R$ to
    be the closed subgroup $\langle \sigma^px_{u-1}\rangle \subset
    N\rtimes \Z_p$. Observe that $R$ is normal since $x_u=1$.
    Finally we form $\G = (N\rtimes\Z_p) /R$ and set $\N$ to be the
    image of $N\rtimes \{1\}$ in $\G$. Since $\N\simeq N$ as
    pro-$p$-$G$ operator groups, we identify them in what follows.

    By construction $\N$ is a maximal closed subgroup of $\G$ which
    is abelian of exponent $p$.   Hence $\G$ is a $T$-group. Since
    the image of $\sigma$ in $\Gamma$ has order $p^2$ we see that
    $\G$ is not of exponent $p$. From the decomposition of $N$, we
    obtain by Lemma~\ref{le:t}\eqref{it:lt4} that the invariants
    $t_i$ are as desired.  It remains only to show that $u$ is as
    given. By Lemma~\ref{le:t}\eqref{it:lt5} we have $\G^p = \langle
    x_{u-1}\rangle\cdot \G_{(p)}$.  From
    Lemma~\ref{le:t}(\ref{it:lt2},\ref{it:lt3}) we calculate that
    $x_{u-1}\in \G_{(u)}$ and $x_{u-1}\not\in \G_{(u+1)}$.  Hence
    $u$ is as desired.

    \emph{The case $u=p$}.  We proceed analogously.  Let $G$ be a
    group of order $p$ and $\Ic_i$, $i=1, 2, \dots, p$ sets
    with cardinalities $\vert \Ic_i\vert$ satisfying
    \begin{align*}
        \vert \Ic_i\vert &= t_i, \qquad 2\le i\le p\\
        \vert \Ic_i\vert &= t_i, \qquad i=1 \text{\ and some\ }
        t_j\neq 0,\ 2\le j\le p\\
        1+\vert \Ic_i\vert &= t_i, \qquad i=1 \text{\ and\ } t_j=0,
        \ 2\le j\le p.
    \end{align*}
    Set $N = M_1^{\Ic_1}\times M_2^{\Ic_2} \times \cdots \times
    M_p^{\Ic_p}$, where $N$ is endowed with the product topology.
    Set $\G=N\rtimes G$, $\N=N\rtimes \{1\}$, and let $\sigma$ be a
    generator of $G$.

    By construction $\N$ is a maximal closed subgroup of $\G$ which
    is abelian of exponent $p$.  Hence $\G$ is a $T$-group and is
    abelian if and only if $t_i=0$ for $i\ge 2$, and in this case
    $\G$ is of exponent $p$.  From the decomposition of $N$ we
    obtain by Lemma~\ref{le:t}\eqref{it:lt4} that the invariants
    $t_i$ are as desired.  It remains only to show that $u$ is as
    given.  By Lemma~\ref{le:t}\eqref{it:lt5} we have $\G^p =
    \G_{(p)}$ and $u=p$, as desired.

    Now we show that $T$-groups are uniquely determined up to
    isomorphism by the invariants $t_i$ and $u$.  Let $T$ be an
    arbitrary $T$-group with invariants $t_i$, $i=1, 2, \dots, p$,
    and $u$, and $N$ be a maximal closed subgroup of $T$ that is
    abelian of exponent dividing $p$.  Let $G=T/N$.  From
    Lemma~\ref{le:t}(\ref{it:lt4}) the structure of $N$ as an $\Fp
    G$-module is determined up to isomorphism, and $T$ is an
    extension of $N$ by $G$.  Let $\sigma\in T\setminus N$.  Since
    $N$ is of index $p$ we have that $\sigma^p\in N$, and since
    $\bar\sigma\in G$ generates $G$ we have that $\sigma^p\in N^G$.
    It remains only to determine the isomorphism class of $N$ as an
    $\Fp G$-module with a distinguished factor $X$ such that
    $\sigma^p\in X^G$.

    \emph{Case 1}: $u = p$, $t_1\ge 1$, and $t_i = 0$ for all $2 \le
    i \le p$. From Lemma~\ref{le:t1}(\ref{it:tc1}) we see that $T$
    is abelian so that $T_{(2)}=\{1\}$.  Then from the definition of
    $u$ we see that $T$ has exponent $p$.  We deduce that $T \simeq
    M_1^{\Ic_1} \times G$ where $G$ is a cyclic group of order $p$
    and $t_1 = \vert \Ic_1 \vert + 1$.  Thus $T$ is determined by
    the invariants.

    \emph{Case 2}: $t_1 \ge 1$, $t_i = 0$ for all $2 \le i \le
    p$, and $u = 1$.  Again we deduce that $T$ is abelian.  From the
    definition of $u$ we see that $T$ does not have exponent $p$,
    and so by the definition of $T$-group, $T$ has exponent $p^2$.  We
    deduce that $N = X \times \tilde N$ where $\sigma^p$ generates
    an $\Fp G$-module $X$ isomorphic to $M_1$ and $\tilde N\simeq
    M_1^{\Ic'}$ with $\vert \Ic' \vert+1 = t_1$. Thus $T$ is
    determined by the invariants.

    \emph{Case 3}: $t_i \ne 0$ for some $i$ with $2 \le i\le p$. By
    Lemma~\ref{le:t1}\eqref{it:tb1} we see that $T$ is nonabelian.
    By Lemma~\ref{le:t}(\ref{it:lt5}), $T^p= \langle\sigma^p\rangle
    \cdot T_{(p)}$.  From Lemma~\ref{le:t}\eqref{it:lt3} and the
    definition of $u$, we deduce that $\sigma^p\in
    N^{(\sigma-1)^{u-1}}$ and, if $u<p$, that $\sigma^p\not\in
    N^{(\sigma-1)^{u}}$.  When $u<p$ we obtain
    \begin{equation*}
        \sigma^p \in \big( N^G\cap N^{(\sigma-1)^{u-1}}\big)
        \setminus \big(N^G \cap N^{(\sigma-1)^u}\big),
    \end{equation*}
    while when $u=p$, we have $\sigma^p \in N^{(\sigma-1)^{p-1}}$.

    If $u<p$ then by Proposition~\ref{pr:tc} we have $t_u\ge 1$, and
    using the same argument as in the proof of
    Lemma~\ref{le:t}\eqref{it:lt6}, we see that $N$ contains a
    distinguished direct factor $X\simeq M_u$ such that $X^G=\langle
    \sigma^p\rangle$. We deduce that $N = X \times \tilde N$ for
    $\tilde N \simeq M_u^{\Ic'} \times \prod_{i\neq u} M_i^{\Ic_i}$
    for sets $\Ic_i$, $i\neq u$, and $\Ic'$ such that
    $\vert\Ic_i\vert = t_i$, $i\neq u$, and $1+\vert \Ic'\vert =
    t_u$.  Then $T$ is determined by its invariants.

    If $u=p$, we claim that without loss of generality we may assume
    that $\sigma^p=1$.  Since $\sigma^p\in N^{(\sigma-1)^{p-1}}$,
    let $\nu\in N$ such that $\sigma^p=\nu^{(\sigma-1)^{p-1}}$ and
    set $\tau = \sigma\nu^{-1}$.  Then $\tau\in T\setminus N$ and
    $\tau^p = \sigma^p (\nu^{-1})^{1+\sigma+\cdots+\sigma^{p-1}}$.
    Since $(\sigma-1)^{p-1}=1+\sigma+\cdots+\sigma^{p-1}$ in $\Fp G$
    we deduce that $\tau^p=1$, as desired.  Hence $T=N\rtimes G$ and
    is determined by its invariants.
\end{proof}

\begin{lemma}\label{le:freet}
    Suppose that $\G$ is a profinite group such that $\G(p)$ is a
    free pro-$p$-group of (possibly infinite) rank $n$, and let $\N$
    be a normal subgroup of $\G$ of index $p$.  Then the invariants
    of the $T$-group $\G/\N^p[\N,\N]$ are given by
    \begin{align*}
        t_i &= \begin{cases} 1, &i=1\\
            0, &2\le i< p\\
            n-1, &i=p, \ n<\infty\\
            n, &i=p, \ n \text{\ an\ infinite\ cardinal}
        \end{cases}\\
        u &= 1.
    \end{align*}
\end{lemma}
\noindent We observe the convention that a free pro-$p$-group is a
free pro-$p$-group on a positive number of generators.

\begin{proof}
    Since $\G/\N^p[\N,\N]=\G(p)/\Phi(\N(p))$, we may assume without
    loss of generality that $\G$ is a free pro-$p$-group. Let
    $\{\sigma_i\}_{i\in \Ic} \cup \{ \sigma \}$ be a minimal
    generating set of $\G$ such that $\sigma \not\in \N$ and all
    $\sigma_i\in \N$.  Write $\bar\sigma_i$ for the image of
    $\sigma_i$ in $T:=\G/\N^{(2)}$.  We use the analogue of the
    Kurosh subgroup theorem in the context of pro-$p$-groups.

    For each $i\in \Ic$, let $Q_i$ be the closed subgroup of $\G$
    generated by $\sigma_i$, and let $Q$ be the closed subgroup
    generated by $\sigma$.  Then $\G$ is the free product
    $\left(\star_{i \in \Ic}\  Q_i \right) \star Q$ in the category
    of pro-$p$-groups.  Applying \cite[Theorem~4.2.1]{NSW} we see
    that
    \begin{equation*}
        \N = \left(\bigstar_{0\le i,j\le
        p-1} \ Q_i^{\sigma^j}\right) \star (Q^p \cap \N),
    \end{equation*}
    where $Q_i^{\sigma^j} = \sigma^{-j} Q_i \sigma^j$ and $Q^p$ is
    the closed subgroup of $Q$ generated by $\sigma^p$.  We deduce
    that $\N/ \N^{(2)}$ is a topological product of copies of
    $\Z/p\Z$ corresponding to the generators of the factors of
    $\N$ above.

    For each $i \in \Ic$ we obtain the $\Fp G$-submodule of $\N$
    \begin{equation*}
        Y_i := \langle \bar\sigma_i,\bar\sigma^{-1} \bar\sigma_i
        \bar\sigma, \dots,\bar\sigma^{-(p-1)} \bar\sigma_i
        \bar\sigma^{p-1} \rangle =
        \langle \bar\sigma_i,\bar\sigma_i^{\bar\sigma-1},\dots,
        \sigma_i^{(\bar\sigma-1)^{p-1}} \rangle\simeq M_p,
    \end{equation*}
    where the equality follows from the fact that $Y_i$ is the
    cyclic $\Fp G$-module generated by $\bar\sigma_i$. We deduce
    that $\N/\N^{(2)} \simeq M_1 \times M_p^{\Ic_p}$, where $\vert
    \Ic_p \vert + 1 = n$.  We calculate that $u=1$ since $\sigma^p
    \in M_1$ and $\sigma^p \notin T_{(2)}$. Finally
    Lemma~\ref{le:t}\eqref{it:lt4} tells us that $t_i$, $i =
    1,2,\dots, p$, have the prescribed values.
\end{proof}

\begin{lemma}\label{le:freegf}
    Let $S$ be a free pro-$p$-group.  Then there exists a field $F$
    of characteristic $0$ such that $G_F\simeq S$.
\end{lemma}

\begin{proof}
    First we let $F_0$ be any algebraically closed field of
    characteristic $0$ with cardinality greater than or equal to $d
    = \dim_{\Fp} H^1(S,\Fp)$; such a field exists as each field
    admits an algebraic closure.  Set $F_1:=F_0(t)$.  By \cite{Do}
    (see also \cite{vDR}), $G_{F_1}$ is a free profinite group, and
    we let $P$ denote a $p$-Sylow subgroup of $G_{F_1}$.  By
    \cite[Corollary~7.7.6]{RZ}, $P$ is a free pro-$p$-group.  Let
    $F_2$ be the fixed field of $P$.

    Consider the set $A$ of linear polynomials $t-c\in F_1$ for
    $c\in F_0$.  Because $F_0[t]$ is a unique factorization domain
    we obtain that the classes of these polynomials in $F_1^\times/
    F_1^{\times p}$ are linearly independent over $\Fp$.  Choose a
    subset of cardinality $d$ of $A$, and let $V$ be the vector
    subspace of $F_1^\times/F_1^{\times p}$ generated by the classes
    of the elements of $A$.  Thus $\dim_{\Fp} V=d$.

    Since $[F_2:F_1]$ is coprime to $p$ (for results on
    supernatural numbers, see \cite[Chapter~1]{SGC}), we deduce
    that $V$ injects into $F_2^{\times}/F_2^{\times p}$.  Let $W$
    denote this image.  Now let $F$ be a maximal algebraic field
    extension of $F_2$ such that $W$ injects into
    $F^\times/F^{\times p}$.

    We claim that the image $i(W)$ of $W$ in $F^{\times}/F^{\times
    p}$ is in fact $F^{\times}/F^{\times p}$.  Suppose not.  Then
    there exists a class $[f]\in F^{\times}/F^{\times p}$, $f\in
    F^\times$, such that $[f]\not\in i(W)$.  Then the subgroup
    $\langle [f]\rangle$ is of order $p$ with trivial intersection
    with $i(W)$.  Set $L=F(\root{p}\of{f})$.  We have that $W$
    injects into $L^\times/L^{\times p}$, but $L$ is strictly larger
    than $F$, a contradiction.  The rank of $G_F$ is then
    $\dim_{\Fp} H^1(G_F,\Fp) = \dim_{\Fp} F^{\times}/F^{\times p} =
    d$.
\end{proof}

\begin{proof}[Proof of Theorem~\ref{th:t}]
    \emph{The case $p=2$}. Let $u\in \{1,2\}$, $t_1$, and $t_2$ be
    invariants of a $T$-group $T$.  By Proposition~\ref{pr:tc}, we
    have that $t_1\ge 1$ if $u=1$, while if $u=2$ then either
    $t_1\ge 1$ or $t_2\ge 1$.

    \emph{Case 1}: $T$ is not of exponent $2$.  By
    Lemma~\ref{le:t1}(\ref{it:tb2}) we find that either $u=1$ or
    $t_2\ge 1$.  Let $G$ be a group of order $2$ and
    $N=M_1^{\Ic_1}\times M_2^{\Ic_2}$ a topological $\Ft G$-module
    for some sets $\Ic_1$ and $\Ic_2$ satisfying
    $\vert\Ic_1\vert=t_1$ and $\vert\Ic_2\vert=t_2$, and let
    $M=N^\vee$.  From \cite[Corollary~2]{MSw2}, there exists a field
    extension $E/F$ of fields of characteristic not $2$ with
    $\Gal(E/F)\simeq G$ such that $H^1(G_E,\Ft)\simeq
    E^\times/E^{\times 2} \simeq M$ as $\Ft G$-modules if and only
    if there exist $\Upsilon\in \{0,1\}$ and cardinalities $d$ and
    $e$ such that $t_1+1=2\Upsilon+d$; $t_2+\Upsilon=e$; if
    $\Upsilon=0$ then $d\ge 1$; and if $\Upsilon=1$ then $e\ge 1$.
    Moreover, the resulting extension $E/F$ will have $-1\in
    N_{E/F}(E^\times)$ if and only if $\Upsilon=1$.  Finally, by
    \cite[proof of Theorem~1]{MSw2}, we may choose $E/F$ such that
    $G_F$ is a pro-$2$-group.

    If $u=1$ then we set $\Upsilon=1$ and $e=t_2+\Upsilon$.  Since
    $t_1\ge 1$ we may choose $d\ge 0$ such that $2\Upsilon+d=t_1+1$.
    We have $e\ge 1$ and hence we have satisfied the conditions for
    a field extension $E/F$ with $G_E/G_E^{(2)}\simeq M^\vee \simeq
    N$. Since $\Upsilon=1$ we have $-1\in N_{E/F}(E^\times)$ and by
    Albert's criterion \cite{A}, $E/F$ embeds in a cyclic extension
    $E'/F$ of degree $4$.  We see that $E'/E$ is degree $2$ and
    hence $\Gal(E'/E)$ is a quotient of $G_E/G_E^{(2)}$.  Let
    $\Gal(E/F)=\langle \bar\sigma\rangle$. Then we have the
    commutative diagram
    \begin{equation*}
        \xymatrix{ 1 \ar[r]& G_E/G_E^{(2)} \ar[r] \ar@{>>}[d] &
        G_F/G_E^{(2)} \ar[r] \ar@{>>}[d] & G \ar[r] \ar[d]^{=} & 1
        \\ 1 \ar[r] & M_1 \ar[r] & H \ar[r] & G \ar[r] & 1}
    \end{equation*}
    in which $\bar\sigma$ lifts to an element of order $4$.  We
    deduce from Lemma~\ref{le:t}\eqref{it:lt6} that the invariant
    $u$ for $T_{E/F}$ is $1$. By Lemma~\ref{le:t1}\eqref{it:tb2},
    $T_{E/F}$ is not of exponent $2$.  From
    Lemma~\ref{le:t}\eqref{it:lt4} and the isomorphism
    $G_E/G_E^{(2)}\simeq N\simeq M_1^{\Ic_1} \times M_2^{\Ic_2}$ we
    find that the invariants $t_1$ and $t_2$ of $T_{E/F}$ match
    those of $T$.  Hence the $T$-groups $T$ and $T_{E/F}$ have the
    same invariants and we conclude by Proposition~\ref{pr:tc} that
    $T\simeq T_{E/F}$.

    If $u=2$ then we determined that $t_2\ge 1$.  We take
    $\Upsilon=0$, $d=t_1+1\ge 1$, and $e=t_2$ and obtain a cyclic
    extension $E/F$ as before. Since $\Upsilon=0$ we have $-1\not\in
    N_{E/F}(E^\times)$ and by Albert's criterion \cite{A}, $E/F$
    does not embed in a cyclic extension $E'/F$ of degree $4$.  Let
    $\Gal(E/F)=\langle \bar\sigma\rangle$. We have that there is no
    commutative diagram
    \begin{equation*}
        \xymatrix{ 1 \ar[r]& G_E/G_E^{(2)} \ar[r] \ar@{>>}[d] &
        G_F/G_E^{(2)} \ar[r] \ar@{>>}[d] & G \ar[r] \ar[d]^{=} & 1
        \\ 1 \ar[r] & M_1 \ar[r] & H \ar[r] & G \ar[r] & 1}
    \end{equation*}
    in which $\bar\sigma$ lifts to an element of order $4$.  We
    deduce from Lemma~\ref{le:t}\eqref{it:lt6} that the invariant
    $u$ of $T_{E/F}$ is $2$. Because $t_2\ge 1$, $G_E/G_{E}^{(2)}$
    contains an $\Ft G$-submodule isomorphic to $M_2$, and
    $\bar\sigma$ does not act trivially on $M_2$.  Hence $T_{E/F}$
    is nonabelian.  By Lemma~\ref{le:t}\eqref{it:lt4} and the
    isomorphism $G_E/G_E^{(2)}\simeq N\simeq M_1^{\Ic_1} \times
    M_2^{\Ic_2}$ we deduce that the invariants $t_1$ and $t_2$ of
    $T_{E/F}$ match those of $T$.  Hence the $T$-groups $T$ and
    $T_{E/F}$ have the same invariants and we conclude by
    Proposition~\ref{pr:tc} that $T\simeq T_{E/F}$.

    \emph{Case 2}: $T$ has exponent $2$. By
    Lemma~\ref{le:t1}(\ref{it:tb2}) we have that $u=2$ and $t_2=0$,
    and therefore by the restrictions for abstract $T$-groups we
    deduce $t_1\ge 1$.  Let $N=M_1^{\Ic_1}$ be a topological $\Ft
    G$-module for a set $\Ic_1$ satisfying $\vert\Ic_1\vert+1=t_1$.
    Choose $\Upsilon=0$, $d=t_1+1$, and $e=t_2=0$.  Then $d\ge 1$
    and there exists a field extension $E/F$ of fields of
    characteristic not $2$ such that $H^1(G_E,\Ft)\simeq M$.
    Moreover $-1\not\in N_{E/F}(E^\times)$ if and only if
    $\Upsilon=0$.  As before we deduce that $u=2$ for $T_{E/F}$, and
    by Lemma~\ref{le:t}\eqref{it:lt4} we find that $t_2$ of
    $T_{E/F}$ is $0$. Hence by Lemma~\ref{le:t1}(\ref{it:tb2}),
    $T_{E/F}$ has exponent $2$.  By Lemma~\ref{le:t}\eqref{it:lt4}
    again we obtain that $t_1$ for $T_{E/F}$ is equal to $t_1$.
    Hence the $T$-groups $T$ and $T_{E/F}$ have the same invariants
    and we conclude by Proposition~\ref{pr:tc} that $T\simeq
    T_{E/F}$.

    We have shown then that for $p=2$, any abstract $T$-group may be
    realized as the $T$-group of $E/F$ for a cyclic extension of
    degree $2$ of fields of characteristic not $2$.  Hence the
    equivalence holds for all $E/F$.

    \emph{The case $p>2$}.  First we characterize those $T$-groups
    occurring as $T_{E/F}$ for fields $F$ such that the maximal
    pro-$p$-quotient $G_F(p)$ of the absolute Galois group $G_F$ is
    free pro-$p$.  Lemma~\ref{le:freet} tells us that for such a
    field $F$ and a cyclic extension $E/F$ of degree $p$, the
    invariants of $T_{E/F}$ must satisfy $t_1=1$, $t_i=0$ for $2\le
    i<p$, and $u=1$, and that the rank of $G_F(p)$ is one more than
    the invariant $t_p$.  Now suppose that $T$ is a $T$-group with
    invariants $t_1=1$, $t_i=0$ for $2\le i<p$, and $u=1$.  By
    Lemma~\ref{le:freegf} there exists a field $F$ such that $G_F$
    is a free pro-$p$-group of rank $t_p+1$. Letting $\N$ be any
    maximal closed subgroup of $G_F$ and $E=\Fix(\N)$,
    Lemma~\ref{le:freet} tells us that the invariants of $T_{E/F}$
    match those of $T$, and then by Proposition~\ref{pr:tc} we have
    that $T\simeq T_{E/F}$. Therefore the $T$-groups which occur as
    $T_{E/F}$ for fields $F$ with free maximal pro-$p$-quotient $G_F(p)$
    are precisely those for which $t_1=1$, $t_i=0$ for $2\le
    i<p$, and $u=1$.

    Now we characterize which of the remaining $T$-groups occur as
    $T_{E/F}$ for cyclic field extensions $E/F$ of degree $p$ for
    $F$ a field either of characteristic $p$ or containing a
    primitive $p$th root of unity.  Because the remaining $T$-groups
    cannot be $T_{E/F}$ for a cyclic field extension $E$ of a field
    $F$ with $G_F(p)$ a free pro-$p$ maximal pro-$p$-quotient
    $G_F(p)$, by Witt's Theorem we have that the characteristic of
    $F$ is not $p$.  Hence we consider only fields $F$ of
    characteristic not $p$ and containing a primitive $p$th root of
    unity.

    For a cyclic extension $E/F$ of degree $p$, consider the group
    \begin{equation*}
        M_{E/F}:=H^1(G_E,\Fp)\simeq E^\times/E^{\times p}.
    \end{equation*}
    This group $M_{E/F}$ is an $\Fp\Gal(E/F)$-module. Now let $G$
    be an abstract group of order $p$.  Since the particular
    isomorphism $G\simeq \Gal(E/F)$ does not alter the isomorphism
    class of $M_{E/F}$ as an $\Fp G$-module, we may consider all
    such modules $\Fp G$-modules. (See \cite{MSw2}.)

    By \cite[Corollary~1]{MSw2}, $M\simeq M_{E/F}$ as $\Fp
    G$-modules for suitable $E/F$ with $G\simeq \Gal(E/F)$ if and
    only if
    \begin{equation*}
        M = M_1^{\Ic_1} \oplus M_2^{\Ic_2} \oplus M_p^{\Ic_p}
    \end{equation*}
    where the cardinalities $j_1=\vert \Ic_1\vert$, $j_2=\vert
    \Ic_2\vert$, and $j_p=\vert\Ic_p\vert$ satisfy the following
    conditions: $j_1+1= 2\Upsilon+d$, $j_2=1-\Upsilon$, and
    $j_p+1=e$ for some cardinalities $d$, $e$ and $\Upsilon\in
    \{0,1\}$ where $d\ge 1$ if $\Upsilon=0$ and $e\ge 1$.  (Recall
    that further, by the Krull-Schmidt-Azumaya theorem, the
    decomposition above into indecomposable $\Fp G$-modules is
    unique up to equivalence.) Moreover, in the given $E/F$ we have
    $\Upsilon=1$ if and only if $\xi_p\in N_{E/F}(E^\times)$.
    Finally, by \cite[proof of Theorem~1]{MSw2}, we may choose $E/F$
    such that $G_F$ is a pro-$p$-group.

    We observe that the constraints on $j_1$, $j_2$, and $j_p$ are
    then $j_1\ge \Upsilon$ and $j_2=1-\Upsilon$.  Now by duality,
    $H^1(G_E,\Fp)^\vee \simeq G_E/G_E^{(2)}$, and since cyclic $\Fp
    G$-modules are self-dual, we may derive conditions on the
    topological $\Fp G$-module $G_E/G_E^{(2)}$ occurring as maximal
    closed subgroups of $T$-groups $G_F/G_E^{(2)}$, as follows. Set
    $G=G_F/G_E$.

    First we relate $\Upsilon$ and the invariant $u$. We claim that
    for any $T$-group $G_F/G_E^{(2)}$, we have $u\le 2$.  Write
    $E=F(\root{p}\of{a})$ for some $a\in F^\times$.  Let $[e]\in
    E^\times/E^{\times p}$ denote the class of $e\in E^\times$. Then
    $X=\langle [\root{p}\of{a}], [\xi_p]\rangle$ is a cyclic $\Fp
    G$-submodule of $M$ and is isomorphic to $M_i$ for some $i \in \{ 1,2
    \}$. By equivariant Kummer theory (see \cite{W}), $L=E(\root{p^2}
    \of{a},\xi_{p^2})$ is a Galois extension of $F$. Moreover,
    $G_F/G_L$ is a homomorphic image of $G_F/G_E^{(2)}$, since $L/E$
    is an elementary abelian extension.
    Then $L/F(\xi_{p^2})$ is Galois with group $\Z/p^2\Z$, and
    for any $\sigma\in G_F\setminus G_E$, the restriction
    $\sigma_L\in \Gal(L/F)$ restricts to a generator. Hence $1\neq
    \sigma^p\in \Gal(L/E)$, and therefore $1\neq \sigma^p\in
    G_E/G_E^{(2)}$.  We have a commutative diagram
    \begin{equation*}
        \xymatrix{ 1 \ar[r]& G_E/G_E^{(2)} \ar[r] \ar@{>>}[d] &
        G_F/G_E^{(2)} \ar[r] \ar@{>>}[d] & G \ar[r] \ar[d]^{=} & 1
        \\ 1 \ar[r] & M_i \ar[r] & G_F/G_L \ar[r] & G \ar[r] & 1}
    \end{equation*}
    and so by Lemma~\ref{le:t}(\ref{it:lt6}) we deduce that $u\le
    i\le 2$.

    Now we claim that $\Upsilon=1$ if and only if $u=1$. We have
    that $\Upsilon=1$ if and only if $\xi_p\in N_{E/F}(E^\times)$.
    By Albert's criterion \cite{A}, $\Upsilon=1$ if and only if
    $E/F$ embeds in a cyclic extension of $F$ of degree $p^2$, if
    and only if there exists a closed normal subgroup $\tilde
    N\subset G_E$ such that $G_F/\tilde N\simeq \Z/p^2\Z$.  Any such
    closed normal subgroup must contain $G_E^{(2)}$.  Hence we
    deduce that $\Upsilon=1$ if and only if there exists a
    commutative diagram
    \begin{equation*}
        \xymatrix{ 1 \ar[r]& G_E/G_E^{(2)} \ar[r] \ar@{>>}[d] &
        G_F/G_E^{(2)} \ar[r] \ar@{>>}[d] & G \ar[r] \ar[d]^{=} & 1
        \\ 1 \ar[r] & M_1 \ar[r] & H \ar[r] & G \ar[r] & 1}
    \end{equation*}
    in which nontrivial elements of $G$ lift to elements of order
    $p^2$.  By Lemma~\ref{le:t}\eqref{it:lt6}, we deduce that
    $\Upsilon=1$ if and only if $u=1$.

    We have therefore shown that $u\le 2$ and $u=2-\Upsilon$.
    Translating the remaining conditions on $j_1$ and $j_2$, we see
    that $j_1\ge 2-u$ and $j_2=u-1$.  Now by
    Lemma~\ref{le:t}\eqref{it:lt4}, $t_2=j_2$, $t_i=0$ for $3\le
    i\le p-1$, and $t_p=j_p$.  Moreover, $t_1=j_1$ if $T$ is not
    abelian of exponent $p$.  By Lemma~\ref{le:t1}\eqref{it:tb2},
    $T$ is of exponent $p$ if and only if $u=p$.  But we have shown
    that $u\le 2<p$, whence $T$ is not of exponent $p$ and we have
    $t_1=j_1$. By Proposition~\ref{pr:tc}, if $u=1$ then $t_1\ge 1$,
    and since $t_2=u-1$ the conditions for applying
    \cite[Corollary~1]{MSw2} are valid. We conclude that a $T$-group
    $T_{E/F}$ with prescribed invariants subject to conditions (1)
    exists.
\end{proof}

\begin{lemma}\label{le:num}
    Let $T$ be a $T$-group with invariants $t_i$, $i=1, 2, \dots,p$,
    and $u$, and let $N$ be a maximal closed subgroup of $T$ which
    is abelian of exponent dividing $p$. Then
    \begin{multline*}
        1 + \dim_{\Fp} H^1(N,\Fp) = \\ \begin{cases}
                    1 + \sum_{i=1}^p i t_i, &u < p \mbox{ or }
                    t_i > 0 \mbox{ for some } 2 \le i \le p. \\
                    t_1 &\mbox{otherwise.}
                  \end{cases}
    \end{multline*}

    Furthermore, if $t_i > 0$ for some $2 \le i \le p$ or if $u = p$ and
    $t_i = 0$ for $2 \le i \le p$, then $Z(T)$ has exponent $p$ and
    $\dim_{\Fp} H^1(Z(T),\Fp) = \sum_{i=1}^p t_i$.
\end{lemma}

\begin{proof}
     Both statements follow from
     Lemma~\ref{le:t1}(\ref{it:tb1},\ref{it:tb2}) and
     Lemma~\ref{le:t}\eqref{it:lt4}.
\end{proof}

\begin{proof}[Proof of Corollary to Theorem~\ref{th:t}]
    Assume that $T\simeq T_{E/F}$ for some $E/F$ of degree $p$.
    Theorem~\ref{th:t} gives us that $u\in \{1,2\}$, so $u<p$. By
    duality $e=\dim_{\Fp} N/N^{(2)}=\log_p \vert N \vert$, and
    $\log_p \vert T \vert = 1+e$ since $N$ is of index $p$.  By
    Lemma~\ref{le:num} we have $e=\sum_{i=1}^p it_i$.  Since $T$ is
    nonabelian, by Lemma~\ref{le:t1}\eqref{it:tb1} we have that some
    $t_i$, $2\le i\le p$, is nonzero.  Moreover, from the
    facts that $N$ is of
    index $p$, $T$ is nonabelian, and $N$ is abelian of exponent
    $p$, we deduce that $Z(T)\subset N$.  Hence $Z(T)$ is an
    $\Fp$-vector space and $\dim_{\Fp} Z(T)$ is well-defined.
    By Lemma~\ref{le:num} we have $f=\sum_{i=1}^p t_i$.

    If $u=1$ then by Theorem~\ref{th:t} we deduce $t_i=0$, $2\le
    i<p$, and therefore $e\equiv f \bmod (p-1)$, a contradiction.
    If $u=2$ then by Theorem~\ref{th:t} we deduce $t_2=1$ and
    $t_i=0$, $3\le i<p$.  Therefore $e\equiv f+1\bmod (p-1)$, a
    contradiction.
\end{proof}

Now we consider to what extent the maximal closed subgroup $N$ in
the definition of a $T$-group is unique.  For $p>2$, let $H$ denote
the Heisenberg group of order $p^3$, and for $p=2$, let $H$ denote
the dihedral group of order $8$.  In both cases we have the
presentation
\begin{equation*}
    H = \langle a,b \ \vert \
    a^p=b^p=[a,b]^p=1, [a,b] \text{\ is\ central\ in\ }H\rangle
\end{equation*}
Observe that $N_a=\langle a,[a,b]\rangle$ and $N_b=\langle
b,[a,b]\rangle$ are distinct maximal subgroups which are abelian of
exponent $p$.  Now let $T=H\times (\Z/p\Z)^\Ic$, endowed with the
product topology, for any index set $\Ic$.  Then $N_a\times
(\Z/p\Z)^\Ic$ and $N_b\times (\Z/p\Z)^\Ic$ are again distinct
maximal closed subgroups which are abelian of exponent $p$.  We show
these groups, together with the elementary abelian $T$-groups of
order greater than $p$, are the only $T$-groups with more than one
maximal closed subgroup which is abelian of exponent dividing $p$.
Outside of these cases, then, the $N$ in the definition of a
$T$-group is a characteristic subgroup.

\begin{proposition}\label{pr:Nchar}
    Let $T$ be a $T$-group.  The following are equivalent:
    \begin{itemize}
        \item $T$ has at least two maximal closed subgroups which
        are abelian of exponent dividing $p$
        \item $T$ is either elementary abelian of order greater
        than $p$, or the direct product of $H$ and an abelian
        group of exponent dividing $p$.
    \end{itemize}
\end{proposition}

\begin{proof}
    Let $T$ be a $T$-group and suppose that $N$ and $M$ are distinct
    closed index $p$ subgroups of exponent dividing $p$. Then $N$
    and $M$ are normal, $S=M\cap N$ is a normal closed subgroup of
    $T$ of index $p^2$, and $M/S$ and $N/S$ are distinct subgroups
    of $T/S$ of order $p$ with trivial intersection. There exists an
    element $\sigma S\in T/S$ contained in neither $M/S$ nor $N/S$.
    The actions of $T/N$ on $N$ and $T/M$ on $M$ are then both given
    by conjugation by $\sigma$.  We then identify
    $T/N$ and $T/M$ as $G=\langle \sigma N\rangle=\langle \sigma
    M\rangle$.

    Consider the image of $M$ under the natural homomorphism $T \to
    G=T/N$.  Since $M$ is not contained in $N$, we see that the
    image of $M$ is $G$. Hence $M$ contains an element $m=\sigma n$
    for some element $n\in N$. Since $M$ is abelian, $m=\sigma n$
    acts trivially on $M$. Hence $G$ acts trivially on $M \cap N$.
    By Lemma~\ref{le:t}\eqref{it:lt3}, $T_{(2)} =
    M^{(\sigma-1)}=N^{(\sigma-1)}$, whence $N^{(\sigma-1)}\subset
    M\cap N$. Because $G$ acts trivially on $N^{(\sigma-1)}$, we
    obtain that $N^{(\sigma-1)^2}=\{0\}$.

    By Lemma~\ref{le:t}\eqref{it:lt4}, $N = M_1^{\Ic_1} \times
    \cdots \times M_p^{\Ic_p}$ for some sets $\Ic_i$, and the
    invariants $t_i$, $i\ge 2$, of $T$ satisfy $t_i=\vert
    \Ic_i\vert$.  If $p>2$, then from $N^{(\sigma-1)^2}=\{0\}$ we
    deduce $t_i=0$, $i\ge 3$.  Therefore $N = M_1^{\Ic_1} \times
    M_2^{\Ic_2}$ for sets $\Ic_1$, $\Ic_2$.

    We now show that the invariant $u$ of $T$ is equal to $p$.
    Assume first that $p>2$.  By the proof of
    Lemma~\ref{le:t}\eqref{it:lt5}, $1 = m^p=(\sigma n)^p =
    \sigma^pv$ for some $v\in T_{(p)}$. By
    Lemma~\ref{le:t}\eqref{it:lt3}, $T_{(p)}\subset M_p^{\Ic_p}$,
    while we have already shown that $t_p=\vert \Ic_p\vert=0$.
    Lemma~\ref{le:t}\eqref{it:lt5} then gives us that $T$ is of
    exponent $p$, and so by Lemma~\ref{le:t1}, $u=p$, as desired.
    Now assume that $p=2$ and $u=1$.    By
    Lemma~\ref{le:t}\eqref{it:lt5}, $T^2=\langle \sigma^2 \rangle
    T_{(2)}$.  From $u=1$ we see that $\sigma^2\not\in T_{(2)}$.
    Again by the proof of Lemma~\ref{le:t}\eqref{it:lt5},
    $m^2=(\sigma n)^2 = \sigma^2v$ for some $v\in T_{(2)}$.  But
    then $m^2\neq 1$, a contradiction. Hence $u=2$ in the case $p=2$
    as well.

    \textit{Case A}: $N\not\subset Z(T)$.  Since $G$ acts trivially
    on $M\cap N$ and $M\cap N$ is an index $p$ subgroup of $N$, we
    obtain that $t_2=\vert \Ic_2 \vert = 1$.  Hence the invariants
    of $T$ satisfy $t_2=1$, $t_i=0$ for $i\ge 3$, and $u=p$. Since
    the $T$-group $U = (M_1^{\Ic_1} \times M_2) \rtimes G =
    M_1^{\Ic_1} \times H$ has the same invariants, by
    Proposition~\ref{pr:tc}, $T\simeq U$.

    \textit{Case B}: $N\subset Z(T)$.  Since $N$ is index $p$, we
    obtain that $T$ is abelian.  Since $N$ is of index $p$, $T$ must
    be either be itself of exponent $p$ or the direct product of
    $\Z/p^2\Z$ and $(\Z/p\Z)^\Ic$ for some index set $\Ic$. In this
    latter case, observe that no element of order $p^2$ can be
    contained in $N$. Hence every element of $N$ lies in
    $W=p\Z/p^2\Z \times (\Z/p\Z)^\Ic$. But $W$ is of index $p$ in
    $T$ and so $N=M=W$, a contradiction. Hence $T$ is abelian of
    exponent $p$, and since $M$ and $N$ are distinct subgroups of
    index $p$, $\vert T\vert>p$.
\end{proof}

Then Theorem~\ref{th:t} and the determination of the invariants $u$
in the proof of the previous proposition give the following

\begin{corollary}
    Suppose $p>2$.  Then each $T_{E/F}$ contains precisely one
    maximal closed subgroup which is elementary abelian.
\end{corollary}

\section{Proof of Theorem~\ref{th:1}}\label{se:illus}

\begin{proof}[Proof of Theorem~\ref{th:1}]
    Suppose $\G$ is a pro-$p$-group with maximal closed subgroup
    $\N$, and let $T=\G/\N^{(2)}$, $N=\N/\N^{(2)}$ and $G =
    \Gamma/\Delta$.  Write $\bar\sigma$ and $\bar\tau$ for the
    images of $\sigma$ and $\tau$, respectively.

    (1).  Since $\N$ surjects onto $N$ we have that
    $\bar\sigma\not\in N$, $\bar\tau\in N$, and
    \begin{align*}
        {}^e[\bar\sigma,\bar \tau] &\not\in {}^{p-1}[\bar\sigma,
        N]\\
        {}^{e+1}[\bar\sigma,\bar\tau] &=1.
    \end{align*}
    By Lemma~\ref{le:t}(\ref{it:tb2},\ref{it:tb3}) we have
    \begin{equation*}
        1\neq {}^e[\bar\sigma,\bar\tau] \in T_{(e+1)} \setminus
        T_{(p-1)}.
    \end{equation*}
    Moreover, since $N$ is abelian of index $p$, we deduce from
    ${}^{e+1}[\bar\sigma,\bar\tau]=1$ that
    ${}^e[\bar\sigma,\bar\tau]\in Z(T)$.  We obtain that some
    invariant $t_{i}$ for $3\le i<p$ is nonzero. Now if $\G$ were an
    absolute Galois group of a field $F$, then setting $E=\Fix(\N)$
    we obtain $T=T_{E/F}$, and then we have a contradiction to
    Theorem~\ref{th:t}.

    Now assume additionally that there exists a closed normal
    subgroup $\Lambda\subset \N$ of $\G$ such that $\G/\Lambda\simeq
    \Z/p^2\Z$ and $e=1$.  Let $\tilde\sigma$ denote the image of
    $\sigma$ in $T$.  We have a commutative diagram
    \begin{equation*}
        \xymatrix{ 1 \ar[r]& N \ar[r] \ar@{>>}[d] & T \ar[r]
        \ar@{>>}[d] & G \ar[r] \ar[d]^{=} & 1 \\ 1 \ar[r] & M_1
        \ar[r] & H \ar[r] & G \ar[r] & 1}
    \end{equation*}
    with an image of $\tilde\sigma$ in $H$ of order $p^2$.  Hence by
    Lemma~\ref{le:t}\eqref{it:lt6}, the invariant $u$ of $T$ is
    equal to $1$.  Then as before we deduce that some invariant
    $t_i$, $2\le i<p$, is nonzero. Again using Theorem~\ref{th:t},
    $\G$ is not the absolute Galois group of a field.

    (2). We translate the conditions into $T$ and $N$ as before,
    obtaining two elements $[\bar\sigma,\bar\tau_i]$, $i=1, 2$,
    which generate distinct subgroups of $Z(T)\cap T_{(2)}$ with
    trivial intersection with $T_{(p)}$.  We deduce that the sum of
    $t_i$, $2\le i<p$, is at least two.  If some $t_i$, $3\le i<p$,
    is at least one, then Theorem~\ref{th:t} tells us that $\G$ is
    not the absolute Galois group of a field.  Otherwise we have
    $t_2\ge 2$, and then again by Theorem~\ref{th:t}, $\G$ is not
    the absolute Galois group of a field.

    (3). We translate the conditions into $T$ and $N$ as before,
    obtaining
    $\bar\sigma\not\in N$, $\bar\tau\in N$, and
    \begin{align*}
        \bar\sigma^p&\in {}^{2}[\bar\sigma,
        N].
    \end{align*}
    By Lemma~\ref{le:t}(\ref{it:tb2},\ref{it:tb3}) we have
    $\bar\sigma^p\subset T_{(3)}$, and from
    Lemma~\ref{le:t}(\ref{it:lt5}) we deduce that $T^p\subset
    T_{(3)}$.  Hence the invariant $u$ of $T$ is greater than or
    equal to $3$. Again by Theorem~\ref{th:t}, $\G$ is not the
    absolute Galois group of a field.
\end{proof}

\begin{proof}[Proof of Corollary to Theorem~\ref{th:1}]
    Let $\N$ be the closed subgroup of $\G$ generated as a normal
    subgroup by $\sigma_1^p$, $\sigma_2$, and $\sigma_i$ for $i\in
    \Ic$, and let $\Lambda$ be the closed subgroup of $\G$ generated
    as a normal subgroup by $\sigma_1^{p^2}$, $\sigma_2$, and
    $\sigma_i$ for $i\in \Ic$.  Set also $T = \Gamma/\Delta^{(2)}$.
    Examining the quotient of $\G$ obtained by trivializing
    $\sigma_2$ and each $\sigma_i$, $i\in \Ic$, we see that $\N$ is
    maximal and $\G/ \Lambda \simeq \Z/p^2\Z$.

    Now let $e=f-1$.  By passing from $\G$ to $T$ using bars for
    denoting images of elements of $\G$ in $T$, we see that
    ${}^e[\bar\sigma_1,\bar\sigma_2] \notin
    {}^{p-1}[\bar\sigma_1,\bar\Delta]$ in $T$. On the other hand,
    ${}^{e+1}[\sigma_1,\sigma_2]={}^f[\sigma_1,\sigma_2]=1$ in $T$.
    By Theorem~\ref{th:1}(1), $\G$ is not an absolute Galois group.
\end{proof}

\section{Analogues of Theorems of Artin-Schreier and
Becker}\label{se:asb}

Recall that by Artin-Schreier there are no elements of order $p$ in
the absolute Galois groups of a field $F$, unless $p=2$ and $F$ is
formally real. Becker proved that the same holds for maximal
pro-$p$-quotients of absolute Galois groups \cite{Be}. In this
section we show that in certain small quotients of absolute Galois
groups, there are no non-central elements of order $p$ unless $p=2$
and the base field $F$ is formally real.

Let $p$ be a prime and $F$ a field containing a primitive $p$th root
of unity.  We define the following fields associated to $F$ (see
\cite{AGKM} and \cite{MSp1}):
\begin{itemize}
    \item $F^{(2)}$: the compositum of all cyclic extensions of $F$
    of degree $p$
    \item $F^{(3)}$: the compositum of all cyclic extensions of
    $F^{(2)}$ of degree $p$ which are Galois over $F$
    \item ${F^{(2)}}^{(2)}$: the compositum of all cyclic extensions
    of $F^{(2)}$ of degree $p$.
\end{itemize}
Then we set $W_F = \Gal(F^{(3)}/F)$ and $V_F = \Gal({F^{(2)}}^{(2)}
/F)$.  Observe that for a cyclic extension $E/F$ of degree $p$,
$T_{E/F} =\Gal(E^{(2)}/F)$.

\begin{theorem}\label{th:2}
    Let $p$ be a prime and $F$ a field containing a primitive $p$th
    root of unity $\xi_p$.  The following are equivalent:
    \begin{enumerate}
        \item\label{it:w2} There exists $\sigma\in V_F\setminus
        \Phi(V_F)$ of order $p$
        \item\label{it:w1} There exists $\sigma\in W_F\setminus
        \Phi(W_F)$ of order $p$
        \item\label{it:w0a} There exists $\sigma\in V_F \setminus
        \Phi(V_F)$ such that for each cyclic $E/F$ of degree $p$
        its image $\sigma_{E/F} \in T_{E/F}$ has order at most $p$.
        \item\label{it:w3} $p=2$ and $F$ is formally real.
    \end{enumerate}
    If these conditions hold, the elements whose square roots are
    fixed by $\sigma$ form an ordering of $F$.
\end{theorem}

\begin{remark*}
    The fact that $W_F$ in particular is sufficiently large to
    encode field-theoretic properties was brought to light in
    \cite{MSp2}, which explained the Witt ring of a field $F$ in terms
    of $W_F$ when $p=2$.
\end{remark*}

If $W_F$ is a finite group, we may say more.  Recall that for a
group $W$ an element in a cohomology ring $H^*(W,\Fp)$ is called
\emph{essential} if its restriction to all proper subgroups is zero.
Recall also that $H^*(W,\Fp)$ is \emph{Cohen-Macaulay} if it is free
and finitely generated over a polynomial subalgebra of $H^*(W,\Fp)$.

\begin{corollary*}
    Suppose that either $p$ is an odd prime or $p=2$ and $F$ is not
    formally real. Suppose also that $F^\times/F^{\times p}$ is
    finite. Then the cohomology ring of $W_F$ with
    $\Fp$-coefficients is Cohen-Macaulay and contains no nonzero
    essential elements.
\end{corollary*}

\begin{proof}
    By \cite[Theorem~2.1]{AK}, the result follows if we show that our
    hypotheses imply that $W_F$ is finite and that every element
    $\sigma$ of order $p$ in $W_F$ is central.

    Since $F^\times/F^{\times p}$ is finite, $F^{(2)}$ is a finite
    Galois extension of $F$. Considering a chain $F_0 \subset F_1
    \subset \dots \subset F_n = F^{(2)}$, where $[F_{i+1}:F_i] = p$
    for $i = 0,\dots,n-1$, and using \cite[Theorems 1 and 2]{LMS},
    we see that $F^{(2)\times}/(F^{(2)\times})^p$ is also finite.
    Thus $F^{(3)}$ is a finite Galois extension and $W_F$ is finite.
    From the construction of $F^{(3)}$ we obtain that the Frattini
    subgroup $\Phi(W_F)$ of $W_F$ is $\Gal (F^{(3)}/F^{(2)})$.
    Recall that $F^{(3)}$ is generated by cyclic extensions
    $N/F^{(2)}$ of degree $p$ such that $N/F$ is Galois.  Since
    $\Gal(F^{(2)}/F)$ is a finite $p$-group acting on the cyclic
    group $\Gal(N/F^{(2)})$ of order $p$, we deduce that
    $\Gal(N/F^{(2)})$ lies in the center of $\Gal(N/F)$.  Therefore
    $\Phi (W_F)$ is a subgroup of the center $Z(W_F)$ of $W_F$.

    Now let $\sigma$ be any element of order $p$ as in Theorem~A.2.
    Then $\sigma \in \Phi(W_F) \subset Z(W_F)$, as required.
\end{proof}

\begin{remark*}
    The case $p=2$ was obtained in \cite[\S 3]{AKM}.
\end{remark*}

\begin{proof}[Proof of Theorem~\ref{th:2}]
    $(1)\implies (2)$ and $(1)\implies (3)$. Observe first that
    $F^{(3)}\subset {F^{(2)}}^{(2)}$.  Hence we have a natural
    surjection $V_F \twoheadrightarrow W_F$. Assume that $\sigma \in
    V_F \setminus \Phi(V_F)$ has order $p$. Then its image in $W_F
    \setminus \Phi(W_F)$ also has order $p$. Moreover, for any
    cyclic extension $E/F$ of degree $p$, the image $\sigma_{E/F}$
    of $\sigma$ in $T_{E/F}$ has order at most $p$, as follows.
    Since $E^{(2)}$ is contained in ${F^{(2)}}^{(2)}$ we see that
    $\sigma_{E/F}^p= 1$ in $T_{E/F}$.

    $(2)\implies (4)$.  Suppose that $\sigma\in W_F\setminus
    \Phi(W_F)$ has order $p$.  As in the proof of
    Theorem~\ref{th:t}, let $a\in F^\times\setminus F^{\times p}$ be
    arbitrary such that $F(\root{p}\of{a})\not\subset F(\xi_{p^2})$.
    Set $K_a:=F(\root{p}\of{a},\xi_{p^2})\subset F^{(2)}$.  Then
    $[K_a:F(\xi_{p^2})]=p$ and $L_a := F (\root{p^2}\of{a},
    \xi_{p^2})$ is a cyclic extension of degree $p^2$ of
    $F(\xi_{p^2})$.  Moreover, $L_a/F$ is a Galois extension of $F$
    and $L_a \subset F^{(3)}$ since $\Gal(L_a/F)$ is a central
    extension of degree $\le p$ of $\Gal(K_a/F)$. Now if
    $\sigma(\root{p}\of{a}) \ne \root{p}\of{a}$ then $\sigma$ has
    order $p^2$ in $W_F$. Hence $\sigma$ fixes all $\root{p}\of{a}$
    for $a\in F^\times$ with $F\neq F(\root{p}\of{a})\not\subset
    F(\xi_{p^2})$.  Since $\sigma\not\in \Phi(W_F)$, there must
    exist a cyclic extension $E\subset F^{(2)}$ of $F$ which is not
    fixed by $\sigma$.  We deduce that $\xi_{p^2}\not\in F$ and so
    $E=F(\xi_{p^2})$.  If $p$ is odd then $L:=F(\xi_{p^3})$ is a
    cyclic extension of $F$ of degree $p^2$ and $\sigma$ restricts
    to a generator of $\Gal(L/F)\simeq \Z/p^2\Z$, again a
    contradiction, whence $p=2$.  Now let $s \in W_F \setminus
    \Phi(W_F)$ be an element of order $2$. By \cite[proof of
    Theorem~2.7]{MSp1}, the elements of $F$ whose square roots are
    fixed by $s$ form an ordering of $F$.  Hence $(2)\implies (4)$.

    $(4)\implies (1)$.  Suppose that $(4)$ holds.  By
    \cite[Satz~3]{Be}, there exists an ordering of $F$ whose square
    roots are fixed by some element $s$ of order $2$ in the maximal
    pro-$2$ quotient $G_F(2)$ of $G_F$.  Then the restriction of $s$
    to ${F^{(2)}}^{(2)}$ is the required element $\sigma \in V_F
    \setminus \Phi(V_F)$ of order $2$.  Hence $(4) \implies (1)$.

    $(3)\implies (2)$.  Let $\sigma \in V_F \setminus \Phi(V_F)$
    such that for each cyclic extension $E/F$ of degree $p$ the
    image $\sigma_{E/F}$ of $\sigma$ in $T_{E/F}$ has order at most
    $p$. Let $E = F(\root{p} \of{a})$ such that $\sigma$ acts
    nontrivially on $\root{p}\of {a}$ and $L_a =
    F(\root{p^2}\of{a},\xi_{p^2})$.  As above, since the restriction
    of $\sigma_{E/F}$ to $L_a \subset E^{(2)}$ is not of order
    $p^2$, we deduce that $p=2$ and $\sqrt{-1}\not\in F$. Hence $F$
    is not quadratically closed.  If $F$ is real closed then there
    exists precisely one extension $E/F$ of degree $2$, namely
    $F^{(2)}$, and $F^{(2)}=F^{(3)}= {F^{(2)}}^{(2)}$, whence
    $W_F=T_{E/F}=V_F$. Otherwise, $F^{(3)}$ is a compositum of
    extensions $L/F$ such that $\Gal(L/F)$ is either a cyclic group
    of order $4$ or a dihedral group of order $8$. (See
    \cite[Corollary~2.18]{MSp2}.) On the other hand, each such $L$
    lies in $E^{(2)}$ for a suitable quadratic extension $E/F$: each
    $L$ may be obtained as a Galois closure of $E(\sqrt{\gamma})$
    for some $[E:F] = 2$ and $\gamma \in E^\times$.  Therefore the
    restrictions $\sigma_{L/F}$ of $\sigma$ to the extensions $L/F$
    have order $\le 2$, and so the restriction of
    $\sigma_{F^{(3)}/F}$ of $\sigma$ to $F^{(3)}/F$ has order $2$.
    Hence $(3)\implies (2)$.
\end{proof}

\section{Families of Pro-$p$-Groups That are Not Absolute Galois
Groups}\label{se:fam}

\begin{theorem}\label{th:fam}
    Let $p>3$ be prime.  There exists a pro-$p$-$\Z_p$ operator
    group $\Omega$ such that no group of the form
    \begin{equation*}
        \G := ((\Omega \star \Sigma)\rtimes \Z_p)/\Ec,
    \end{equation*}
    where
    \begin{itemize}
        \item $\Sigma$ is any pro-$p$-group with trivial
        $\Z_p$-action, and
        \item $\Ec$ is any normal closed subgroup of
        $(\Omega\star\Sigma)\rtimes \Z_p$ such that
        \begin{equation*}
            \Ec \subset ((\Omega\star\Sigma)\rtimes p\Z_p)^{(3)},
        \end{equation*}
    \end{itemize}
    is an absolute Galois group.
\end{theorem}
\noindent Here $\star$ denotes the free product in the category of
pro-$p$-groups. (Recall that $R^{(n)}$ denotes the $n$th term of the
$p$-descending series of a pro-$p$-group $R$; see
section~\ref{se:t}.)

The proof of Theorem~\ref{th:fam} relies on \cite[Theorem 1]{LMS}
and the existence of pro-$p$-$\Z_p$ operator groups $\Omega$ with
certain cohomological properties, established in the following
proposition. We prove this proposition in sections~\ref{se:fam2} and
\ref{se:fam3}. Recall that for a pro-$p$-group $\G$, the
\emph{decomposable subgroup} of $H^2(\G,\Fp)$ is defined to be the
subgroup generated by cup products of elements of $H^1(\G,\Fp)$:
\begin{equation*}
    H^2(\G,\Fp)^\dec = H^1(\G,\Fp).H^1(\G,\Fp).
\end{equation*}
We say that $H^2(\G,\Fp)$ is \emph{decomposable} if $H^2(\G,\Fp)=
H^2(\G,\Fp)^\dec$.
\begin{proposition}\label{pr:const}
    Let $p>3$ be prime and $C$ be a cyclic group of order $p$.
    \begin{enumerate}
        \item\label{it:const1} There exists a torsion-free
        pro-$p$-$C$ operator group $\Omega_1$ such that as $\Fp
        C$-modules, $H^1(\Omega_1,\Fp) \simeq M_p$ and
        \begin{align*}
            H^2(\Omega_1,\Fp) = H^2(\Omega_1,\F_p)^\dec
            &\simeq M_{p-1}\oplus
            \bigoplus_{i=1}^{(p-3)/2} M_p.
        \end{align*}
        \item\label{it:const2} There exists a torsion-free
        pro-$p$-$C$ operator group $\Omega_2$ such that as $\Fp
        C$-modules, $H^1(\Omega_2,\Fp) \simeq M_2\oplus M_p$ and
        \begin{align*}
            H^2(\Omega_2,\Fp) = H^2(\Omega_2,\F_p)^\dec
            &\simeq M_1\oplus M_3\oplus
            \bigoplus_{i=1}^{(p+1)/2} M_p.
        \end{align*}
    \end{enumerate}
\end{proposition}
\noindent Here we use $M_i$ to denote the unique cyclic $\Fp
C$-module of dimension $i$, and here and in the remaining sections
we write the action of $\Fp C$ multiplicatively instead of
exponentially.

In what follows we freely use basic properties of the transgression
map found in \cite[Chapters~2 and 3]{NSW}, and we abbreviate this
map by $\tra$. We denote by $\inf$ and $\res$ the inflation and
restriction maps of group cohomology, respectively.  Given a free
pro-$p$-group $V$ and a pro-$p$-group $\N$, we say that a surjection
$V\twoheadrightarrow \N$ is a \emph{minimal presentation} of $\N$ if
$\inf:H^1(\N,\Fp)\to H^1(V,\Fp)$ is an isomorphism.

\begin{theorem}\label{th:a2}
    Let $\Delta$ be a pro-$p$-group, $V$ a free pro-$p$-group, and
    $1 \to R \to V \to \N \to 1$ a minimal presentation of $\N$.
    Then we have natural maps
    \begin{enumerate}
        \item\label{it:n1} $H^2(\N,\Fp)\simeq (R/(R^p[R,V]))^\vee$
        \item\label{it:n2} $H^2(\N,\Fp)^\dec \simeq (R/(R\cap
        V^{(3)}))^\vee$.
    \end{enumerate}
\end{theorem}

\begin{proof}
    Set $\Delta^{[2]} := {\Delta}/{\Delta^{(2)}}$ and $V^{[3]}:=
    V/V^{(3)}$.  Since the presentation is minimal,
    $\Delta^{[2]}\simeq{V}/{V^{(2)}}$. We then have the following
    commutative diagram
    \begin{equation*}
        \xymatrix{1 \ar[r] & R \ar[r] \ar[d] & V \ar[r]
        \ar[d] & \Delta \ar[d] \ar[r] & 1 \\ 1 \ar[r] &
        \frac{V^{(2)}}{V^{(3)}} \ar[r] & V^{[3]} \ar[r] &
        \Delta^{[2]} \ar[r] & 1}
    \end{equation*}
    Here the middle and right-hand vertical maps are the natural
    projections. Since the presentation is minimal, $R\subset
    V^{(2)}$ and the left-hand vertical map is the natural factor
    map.

    Passing to the natural maps induced by the Hochschild-Serre-Lyndon
    spectral sequence with coefficients in $\Fp$, we obtain the
    commutative diagram
    \begin{equation*}\Small
        {\xymatrix{H^1(\Delta^{[2]}) \ar[r]^{\inf} \ar[d] &
        H^1(V^{[3]}) \ar[r]^-{\res} \ar[d] &
        H^1\left(\frac{V^{(2)}}{V^{(3)}}\right)^{\N^{[2]}}
        \ar[r]^-{\tra} \ar[d] & H^2(\Delta^{[2]}) \ar[r] \ar[d]
        & \dots \\ H^1(\N) \ar[r]^{\inf} & H^1(V)
        \ar[r]^{\res} & H^1(R)^\Delta \ar[r]^{\tra} &
        H^2(\Delta) \ar[r] & 0.}}
    \end{equation*}
    (Recall that $H^2(V,\Fp)=0$ since $V$ is a free pro-$p$-group.)
    Since the extension
    \begin{equation*}
        1 \to V^{(2)}/V^{(3)} \to V^{[3]} \to \N^{(2)} \to 1
    \end{equation*}
    is central, we deduce that $H^1(V^{(2)}/
    V^{(3)},\Fp)^{\Delta^{[2]}} = H^1(V^{(2)}/ V^{(3)}, \Fp)$.
    Additionally using the fact that the inflation map on the first
    cohomology group in each row of the previous diagram is an
    isomorphism, we may extract the following commutative square,
    with the rightmost transgression map an isomorphism:
    \begin{equation*}
        \xymatrix{H^1\left(\frac{V^{(2)}}{V^{(3)}},\Fp\right)
        \ar[d]^{\tra} \ar[r] & H^1(R,\Fp)^\N \ar[d]^{\tra}_{\simeq} \\
        H^2(\Delta^{[2]},\Fp) \ar[r] & H^2(\N,\Fp) }
    \end{equation*}
    Since $H^1(R,\Fp)^\N \simeq (R/R^p[R,V])^\vee$, we have
    \eqref{it:n1}.

    The leftmost transgression map, however, is also an isomorphism
    by \cite[1.1 and proof]{Ho}.  Now consider the natural map
    $R/R^p[R,V]\to V^{(2)}/V^{(3)}$ of abelian topological groups of
    exponent $p$.  The image of the natural map
    \begin{equation*}
        H^1(V^{(2)}/V^{(3)},\Fp)\to H^1(R,\Fp)^{\N}
    \end{equation*}
    is $H^1(RV^{(3)}/V^{(3)},\Fp)$. We then factor the horizontal
    maps of the commutative square into homomorphisms followed by
    injections:
    \begin{equation*}
        \xymatrix{H^1\left(\frac{V^{(2)}}{V^{(3)}},\Fp\right)
        \ar[d]^{\tra}_{\simeq} \ar[r] &
        H^1(\frac{RV^{(3)}}{V^{(3)}},\Fp) \ar@{^{(}->}[r] &
        H^1(R,\Fp)^\N \ar[d]^{\tra}_{\simeq} \\
        H^2(\Delta^{[2]},\Fp) \ar[r] & H^2(\N,\Fp)^\dec
        \ar@{^{(}->}[r] & H^2(\N,\Fp). }
    \end{equation*}
    We obtain isomorphisms
    \begin{equation*}
        H^1(R/(R\cap V^{(3)}),\Fp) \simeq H^1(RV^{(3)}/V^{(3)},\Fp)
        \simeq H^2(\N,\Fp)^\dec.
    \end{equation*}
    Hence $H^2(\N,\Fp)^\dec \simeq (R/(R\cap V^{(3)}))^\vee$, and
    we have proved \eqref{it:n2}.
\end{proof}

\begin{remark*}
    The fact that $R^p[R,V] \simeq R \cap V^{(3)}$ if and only if $
    H^2(\N,\Fp)$ is decomposable was proved earlier in the case
    $p=2$ in  \cite[Theorem~2]{GM} and \cite[\S 5]{MSp2}.
\end{remark*}

\begin{theorem}\label{th:gammah}
    Let $\G$ and $H$ be pro-$p$-groups with maximal closed
    subgroups $\N$ and $N$, respectively, and $\alpha:\G
    \to H$ a surjection with $\alpha(\N) = N$ and $\Ker
    \alpha \subset \N^{(3)}$.  Write $G$ for $\G/
    \N \simeq H/N$.  Then as $\Fp G$-modules
    \begin{equation*}
        H^2(\N,\Fp)^\dec \simeq H^2(N,\Fp)^\dec.
    \end{equation*}
\end{theorem}

\begin{proof}
    We shall prove that $H^2(\N,\Fp)^\dec \simeq H^2(N,\Fp)^\dec$
    under a natural isomorphism, and from the proof it will follow
    that the isomorphism is $\Fp G$-equivariant.  Let
    \begin{equation*}
        1 \to R \to V \xrightarrow{\theta} \N \to 1
    \end{equation*}
    be a minimal presentation of $\N$ with $V$ a free pro-$p$-group.
    Since $\ker \alpha\subset \N^{(3)}$ and
    $\theta(V^{(3)})=\N^{(3)}$, we may choose a section $W\subset
    V^{(3)}$ of $\ker \alpha$ under the surjection $\theta$.  We
    then obtain a presentation of $N = \alpha(\N)$ as follows:
    \begin{equation*}
        1 \to RW \to V \xrightarrow{\psi} N \to 1
    \end{equation*}
    where $\psi = \alpha \circ \theta$.  Since $RW \subset V^{(2)}$,
    this presentation is minimal.  By
    Theorem~\ref{th:a2}\eqref{it:n2} we have natural isomorphisms
    \begin{align*}
        H^2(\N,\Fp)^\dec &\simeq (RV^{(3)}/V^{(3)})^\vee \simeq
        (RWV^{(3)}/{V^{(3)}})^\vee \simeq H^2(N,\Fp)^\dec.
   \end{align*}
\end{proof}

\begin{lemma}\label{le:lem}
    Let $p$ be prime and $\G$ a nonfree pro-$p$-group which is the
    absolute Galois group of a field $F$.  Then $\chr F\neq p$ and
    a primitive $p$th root of unity $\xi_p$ lies in $F$.
\end{lemma}

\begin{proof}
    Since $\G$ is not free, then by Witt's Theorem (see
    \cite[Theorem~9.1]{Koc} and \cite[pages~207-210]{JLY}), $\chr
    F\neq p$.  If $F$ does not contain a primitive $p$th root of
    unity, then $F(\xi_p)\subset F^s$ and $\Gal(F(\xi_p)/F)$ is a
    quotient of the absolute Galois group $\G$, while $([F(\xi_p):
    F],p)=1$, a contradiction.  Hence $\xi_p\in F$.
\end{proof}

\begin{theorem}\label{th:deltapgroup}
    Let $p>3$ be prime.  Suppose that $\G$ is a pro-$p$-group and
    $\N$ is a maximal closed subgroup of $\G$.

    If the $\Fp(\G/\N)$-module $H^2(\N,\Fp)^\dec$ contains a cyclic
    summand of dimension $i$ with $3\le i<p$, then $\G$ is not an
    absolute Galois group.

    Moreover, if $\G$ contains a normal closed subgroup
    $\Lambda\subset \N$ with $\G/\Lambda\simeq \Z/p^2$, then we may
    take $2\le i<p$ in the same statement.
\end{theorem}

\begin{proof}
    Suppose that $\G$ and $\N$ are pro-$p$-groups satisfying the
    hypotheses. Contrary to our statement, assume that $\G=G_F$ for
    some field $F$.  Since $\N$ is of index $p$, we obtain that $\N
    = G_E$ for some Galois extension $E/F$ of degree $p$.  Because
    $H^2(\N,\Fp)$ is nonzero, $G_E$ is not a free pro-$p$-group and
    therefore neither is $G_F$ \cite[Corollary 3, \S I.4.2]{SGC}.
    Lemma~\ref{le:lem} tells us that $\chr F \ne p$ and that a
    primitive $p$th-root of unity $\xi_p$ lies in $F$.  By
    \cite[Theorem~11.5]{MeSu}, we obtain that $H^2(\N,\Fp)$ is
    decomposable.  Therefore by \cite[Theorem 1]{LMS},
    $H^2(\N,\Fp)^\dec $ contains no cyclic $\Fp(\G/\N)$-summand of
    dimension $i$ with $3\le i<p$, a contradiction.  Moreover, if
    $\Lambda\subset N$ is a normal closed subgroup and
    $\G/\Lambda\simeq \Z/p^2\Z$, then by Albert's criterion \cite{A},
    $\xi_p\in N_{E/F}(E)$.  Letting $E=F(\root{p}\of{a})$ for some
    $a\in F^\times$, we obtain $(a).(\xi_p)=0$ as an element of
    $H^2(\G,\Fp)$.  Then by \cite[Theorem 1]{LMS}, $H^2(\N,\Fp)^\dec
    $ contains no cyclic $\Fp(\G/\N)$-summand of dimension $2$,
    again a contradiction.
\end{proof}

\begin{corollary}\label{co:deltan}
    Let $p>3$ be a prime.  Suppose that $\G$ and $H$ are
    pro-$p$-groups with respective maximal subgroups $\N$ and $N$,
    and $\alpha:\G \to H$ a surjection with $\alpha(\N) = N$ and
    $\Ker \alpha \subset \N^{(3)}$.  Write $G$ for $\G/\N \simeq
    H/N$. If either $H^2 (\N,\Fp)^\dec$ or $H^2(N,\Fp)^\dec$
    contains a cyclic $\Fp G$-summand of dimension $i$ with $3\le
    i<p$, then neither $\N$ nor $H$ is an absolute Galois group.

    Moreover, suppose additionally that $H$ contains a normal closed
    subgroup $\Lambda\subset N$ with $H/\Lambda\simeq \Z/p^2\Z$.
    Then if either $H^2(\N,\Fp)^\dec$ or $H^2(N,\Fp)^\dec$ contains
    a cyclic $\Fp G$-summand of dimension $2$, then $H$ is not an
    absolute Galois group.
\end{corollary}

\begin{proof}
    By Theorem~\ref{th:gammah} we have that $H^2(\N, \Fp)^\dec
    \simeq H^2(N,\Fp)^\dec$ as $\Fp G$-modules. The remainder
    follows from Theorem~\ref{th:deltapgroup}.
\end{proof}

\begin{proof}[Proof of Theorem~\ref{th:fam}]
    Let $\Omega$ be either of the groups $\Omega_i$ described in
    Proposition~\ref{pr:const} and let $\Sigma$ and $\G$ be defined
    as in the statement of the theorem.

    Observe that $\N :=((\Omega\star\Sigma)\times p\Z_p)/\Ec$ is a
    maximal closed subgroup of $\G$ of index $p$.  Let $G=\G/\N$ and
    note that the actions of $G$ and $C$ on $\N$ are identical. By
    Corollary~\ref{co:deltan} it is enough to show that
    $H^2(\N,\Fp)^\dec$ for $\Ec=1$ contains a cyclic $\Fp G$-summand
    $M_i$ with $3\le i<p$. Assume then that $\Ec=1$.

    By \cite[Theorem~4.1.4]{NSW} we calculate $\Fp G$-module
    isomorphisms
    \begin{align*}
        H^1(\N,\Fp) &\simeq H^1(\Omega,\Fp) \oplus H^1(\Sigma,\Fp)
        \oplus H^1(p\Z_p,\Fp) \\
        H^2(\N,\Fp) &\simeq H^2(\Omega,\Fp) \oplus H^2(\Sigma,\Fp)
        \oplus H^1(\Omega,\Fp)\oplus H^1(\Sigma,\F_p).
    \end{align*}
    (Here we use the fact that $H^2(p\Z_p,\Fp)=\{0\}$ and that
    $H^2(A\times B,\Fp)=H^2(A,\Fp)\oplus H^2(B,\Fp)\oplus
    H^1(A,\Fp)\otimes H^1(B,\Fp)$.)

    By Proposition~\ref{pr:const}, $H^2(\Omega,\Fp)$ is
    decomposable. We deduce that $H^2(\Omega,\Fp)$ is a direct
    summand of $H^2(\N,\Fp)^\dec$. From Proposition~\ref{pr:const},
    we obtain that $H^2(\Omega,\Fp)$ contains an $\Fp G$-summand
    $M_i$ with $3\le i\le p-1$ and we are done.
\end{proof}

The foregoing proof excludes the families of groups $\G$ from the
class of absolute Galois groups by using the $n=2$ case of
\cite[Theorem 1]{LMS}.  When $\Omega=\Omega_1$, the $n=1$ case
cannot be used directly to exclude the groups $\G$. Moreover,
neither can Theorem~\ref{th:t} be used to exclude the groups $\G$,
as follows. Setting $T=\G/\N^{(2)}$ and using duality (see
\cite[Lemma~2.9.4 and Theorem 2.9.6]{RZ}), we obtain from
Lemma~\ref{le:t}\eqref{it:lt4} the invariants $t_p=1$, $t_i=0$ for
$2\le i<p$, and $t_1=1+\dim_{\Fp} H^1(\Sigma,\Fp)$. Observing that
$\Lambda:=((\Omega\star \Sigma)\times p^2\Z_p)/\Ec$ is a normal
subgroup of index $p^2$ with quotient $\Z/p^2\Z$ and contained in
$\N$, we deduce that $u=1$.  But this set of invariants is permitted
by Theorem~\ref{th:t}.

When $\Omega=\Omega_2$, the situation is different.  We see that
$H^1(\N,\Fp)$ does not contain any indecomposable $\Fp C$-modules of
dimension $i$ for $3\le i<p$ and so the groups $\G$ cannot be
excluded by a simple observation of impermissible indecomposable
$\Fp C$-summands in the $n=1$ case of \cite[Theorem 1]{LMS}.
However, since $u=1$ as in the case of $\Omega=\Omega_1$, we may
deduce that the invariant $\Upsilon_2$ of the $n=1$ case of
\cite[Theorem 1]{LMS} is zero, as follows.  Suppose $\G$ is the
absolute Galois group of a field $F$.  Since $\G$ is a pro-$p$-group
which is not free, by Lemma~\ref{le:lem} we deduce that $F$ is of
characteristic not $p$, containing a primitive $p$th root of unity
$\xi_p$.  Hence $E:=\Fix(\N)=F(\root{p}\of{a})$ for some $a\in
F^\times$.  But since $E/F$ embeds in a cyclic extension of degree 4
given by $\Fix(\Lambda)$, we have that $\xi_p\in N_{E/F}(E^\times)$
by \cite{A} and so $(a).(\xi_p)=0$, whence $\Upsilon_2=0$. In a
similar manner, $T=\G/\N^{(2)}$ violates the conditions of
Theorem~\ref{th:t} since $t_2\ge 1$ when $u=1$.

\begin{remark*}
    We observe that the fact that $\Omega_i\rtimes \Z_p$, $i=1, 2$,
    is not an absolute Galois group could also be deduced from the
    main results in \cite{K}, using different methods. However, the
    fact that each $((\Omega_i\times \Sigma)\rtimes \Z_p)/\Ec$ is
    not an absolute Galois group does not follow from \cite{K}.
\end{remark*}

\section{The Group $\Omega_1$ and its Cohomology}\label{se:fam2}

In this section we prove
Proposition~\ref{pr:const}\eqref{it:const1}.

Let $D=\langle g \mid g^p=1 \rangle$ be a cyclic group of order $p$,
and let $\Z_p D$ be the $p$-adic group ring.  Because we will form
group extensions of the additive group of $\Z_p G$, to avoid
confusion we write the group multiplicatively and denote it by
$\GG$, where the element $\bar g_i$ of $\GG$ corresponds to the
element $g^i$ of $\Z_p D$. (We reserve the symbol $g_i$ for elements
of the group extension $\Omega_1$, which we define later.) Thus
$\bar g_0, \bar g_1, \dots, \bar g_{p-1}$ form a topological
generating set of $\GG$. We shall always read the suffixes mod $p$,
so that for example if $i=p-1$ then $\bar g_{i+1}$ denotes $\bar
g_0$.

Now let $C=\langle \sigma \mid \sigma^p=1\rangle$ be another group
of order $p$, acting on $\GG$ via $\sigma(\bar g_i)=\bar g_{i-1}$.
In this way $\GG$ and $H^1(\GG,\Z_p)$ are free $\Z_p C$-modules, and
$H^1(\GG,\Z_p)$ has a topological generating set $y_0, y_1, \dots,
y_{p-1}$ dual to $\bar g_0, \bar g_1, \dots, \bar g_{p-1}$.  Observe
that $\sigma(y_i)=y_{i+1}$.

Next let $\H=\Z_p$ be a trivial $\Z_p C$-module with generator $h$,
and let $z\in H^1(\H,\Z_p)$ be dual to $h$.  We define a nonsplit
central extension $\Omega_1$ of $\H$ by $\GG$ as follows. The group
$H^2(\GG,\H) = \bigwedge^2 H^1(\GG,\H)$ is a free $\Z_p C$-module of
rank $p(p-1)/2$ with free generators $y_0y_j$ for $1\le j\le
(p-1)/2$. Consider the element
\begin{equation*}
    y := (1+\sigma+\dots+\sigma^{p-1})y_0y_1 = y_0y_1+y_1y_2+\dots
    +y_{p-1}y_0.
\end{equation*}
We let $\Omega_1$ be the central extension of $\H$ by $\GG$
corresponding to $y\in H^2(\GG,\H)$:
\begin{equation*}
    1\to \H \to \Omega_1 \to \GG\to 1.
\end{equation*}

The group $\Omega_1$ is a torsion-free nilpotent pro-$p$-group of
Hirsch length $p+1$. Examining the standard correspondence of group
extensions with $H^2$ gives us that for suitable representatives
$g_i\in \Omega_1$ with respective images $\bar g_i\in \GG$, we have
the relations
\begin{align*}
    [g_i,g_j] &= h, \quad j=i+1\\
    [g_i,g_j] &= 1, \quad j\neq i+1, \ i\neq j+1\\
    [g_i,h] &= 1, \quad \text{all\ }i.
\end{align*}
Observe that the action of $C$ on $\GG$ may be extended to
$\Omega_1$ by $\sigma(g_i)=g_{i-1}$, $0\le i< p$, and $\sigma(h)=h$.

The $E_2$ page of the spectral sequence
\begin{equation*}
    H^s(\GG,H^t(\H,\Z_p)) \Rightarrow H^{s+t}(\Omega_1,\Z_p)
\end{equation*}
is generated by the anticommuting elements $y_i\in E_2^{1,0}$ and
$z\in E_2^{0,1}$, and the differential $d_2$ is given by
$d_2^{0,1}(z)=y$.  (See \cite[Theorem~2.1.8]{NSW}.) The same holds
for the $E_2$ page of the spectral sequence with coefficients in
$\Fp$
\begin{equation*}
    H^s(\GG,H^t(\H,\Fp)) \Rightarrow H^{s+t}(\Omega_1,\Fp),
\end{equation*}
and we write $\bar y_i$ and $\bar z$ for the reduction mod $p$ of
the classes $y_i$ and $z$. The picture below indicates the bases of
relevant cohomology groups on this $E_2$ page.
\[ \begin{array}{|c|c|c|c}
 & & & \\ \hline
\bar z & \bar z.\bar y_i & \bar z.\bar y_i.\bar y_j & \\ \hline 1 &
\bar y_i & \bar y_i.\bar y_j & \bar y_i.\bar y_j.\bar y_k \\ \hline
\end{array} \]

First, we observe that $d_2^{0,1}$ is injective on
$E^{0,1}=H^1(\H,\F_p)$, so that by the Five-Term Exact Sequence
\cite[Proposition~2.1.2]{NSW}, $H^1(\Omega_1,\Fp)\simeq
H^1(\GG,\Fp)$ as $\Fp C$-modules. We obtain then that
$H^1(\Omega_1,\Fp)\simeq M_p$.

Now we turn to the calculation of $H^2(\Omega_1,\Fp)$.  We claim
first that $d_2^{1,1}$ is injective, as follows.  It is enough to
show that $d_2^{1,1}(\bar z.\bar x)=0$ implies that $\bar x=0$ for
$\bar x\in H^1(\GG,\Fp)$.  Write $\bar x=\sum_{i=0}^{p-1} c_i \bar
y_i$, $c_i\in \Fp$.  Since $p>3$, the set of elements $\bar y_i.\bar
y_{i+2}.\bar y_{i+3}$, $0\le i<p$, is $\Fp$-independent and may be
expanded to an $\Fp$-basis of $E_2^{3,0}=H^3(\GG,\Fp)$ consisting of
threefold products of distinct $\bar y_i$. Consider the coefficient
$\beta_i$ of $\bar y_i.\bar y_{i+2}.\bar y_{i+3}$ in an expansion of
$d_2^{1,1}(\bar z.\bar x)$. Since $p>3$, the only consecutive pair
of indices in $\bar y_i.\bar y_{i+2}.\bar y_{i+3}$ is $\{i+2,i+3\}$.
Such a consecutive pair appears the product of the term $\bar
y_{i+2}.\bar y_{i+3}$ in $\bar y$ with the term $c_i\bar y_i$ of
$\bar x$.  All other nonzero products of terms of $\bar y$ with
terms of $\bar x$ that contain a consecutive pair then contain more
than one consecutive pair.  Hence $\beta_i$ depends only on $c_i$.
If all $\beta_i=0$, then each $c_i=0$ and therefore $d_2^{1,1}$ is
injective.

Since $d_2^{1,1}$ is injective, we obtain that $E_3^{1,1} =
E_\infty^{1,1}=0$.  Since $E_2^{0,2}=0$ we have $E_\infty^{0,2}=0$.
Next observe that $\bar y = (\sigma-1)^{p-1}\bar y_0.\bar y_1$ and
so $\bar y\in H^2(\GG,\Fp)^C$.  Since $E^{2,0}=H^2(\GG,\Fp)$ is the
direct sum of free $\Fp C$-modules $M(0,i)\simeq M_p$ on generators
$\bar y_0.\bar y_i$, $1\le i\le (p-1)/2$, we deduce that
\begin{equation*}
    E_\infty^{2,0}=E_3^{2,0}=(M(0,1)/M(0,1)^C) \oplus
    \bigoplus_{i=2}^{(p-1)/2} M(0,i).
\end{equation*}
We conclude that $H^2(\Omega_1,\Fp) \simeq M_{p-1} \oplus
\frac{(p-3)}{2} M_p$.

Finally, since $H^1(\Omega_1,\Fp)\simeq H^1(\GG,\Fp)$ and
$H^2(\Omega_1,\Fp)$ is a quotient of $E_2^{2,0}=H^2(\GG,\Fp)=
\bigwedge^2 H^1(\GG,\Fp)$, which is decomposable, we deduce that
$H^2(\Omega_1,\Fp)$ is decomposable as well.

\section{The Group $\Omega_2$ and its Cohomology}\label{se:fam3}

In this section we prove
Proposition~\ref{pr:const}\eqref{it:const2}.

As before let $D=\langle g \mid g^p=1 \rangle$ be a cyclic group of
order $p$, $\GG=\Z_p D$ the $p$-adic group ring with topological
generating set $\bar g_i:=g^i$ for $0\le i<p$, $C=\langle \sigma
\mid \sigma^p=1\rangle$ another group of order $p$ acting on $\GG$
via $\sigma(\bar g_i)=\bar g_{i-1}$, and $H^1(\GG, \Z_p)$ the
resulting free $\Z_p C$-module with topological generating set $y_0,
y_1, \dots, y_{p-1}$ dual to $\bar g_0, \bar g_1, \dots, \bar
g_{p-1}$. Again $\sigma(y_i)=y_{i+1}$.

Now let $H=\langle h \mid h^p=1 \rangle$ be another cyclic group of
order $p$ and $\Z_p H$ be its $p$-adic group ring. As with $\Z_p D$,
we write $\Z_p H$ multiplicatively, denoting $h^i$ by $h_i$. We then
set $\H$ to be the subgroup of $\Z_p H$ generated by the elements
$h_ih_j^{-1}$.  Then $\H$ and $H^1(\H,\Z_p)$ are free $\Z_p$-modules
of rank $p-1$.  Let $z_0,z_1,\dots, z_{p-1}\in H^1(\H,\Z_p)$ be the
images under the surjection $H^1(\Z_p H,\Z_p)\twoheadrightarrow
H^1(\H,\Z_p)$ of the topological generating set of $H^1(\Z_p
H,\Z_p)$ dual to $h_0, h_1, \dots, h_{p-1}\in \Z_p H$.  Since $\H$
is generated by elements of the form $h_ih_j^{-1}$, the elements
$z_i$ satisfy the relation $z_0 + z_1 + \dots + z_{p-1}=0$.
Moreover, the set $\{h_0h_{p-1}^{-1}, h_1h_{p-1}^{-1}, \dots,
h_{p-2}h_{p-1}^{-1}\}$ is dual to the set $\{z_0, z_1, \dots,
z_{p-2}\}$.

We define a nonsplit extension $\Omega_2$ of $\H$ by $\GG$ as
follows.  Writing $\H$ as a free $\Z_p$-module with $\Z_p$-base
$h_ih_{p-1}^{-1}$, $0\le i<p-1$, we have
\begin{equation*}
    H^2(\GG,\H)=\bigoplus_{i=0}^{p-2} H^2\big(\GG,\langle
    h_ih_{p-1}^{-1}\rangle\big)\simeq \bigoplus_{i=0}^{p-2}
    \left(\bigwedge^2 H^1(\GG,\Z_p)\right).
\end{equation*}
With these identifications, consider the element $y$ defined as
\begin{multline*}
    (-y_{0}.y_{1}+3y_{1}.y_2-3y_2.y_3+y_3.y_4,
    -y_{1}.y_2+3y_2.y_3-3y_3.y_4+y_4.y_5, \\ \dots,
    -y_{p-2}.y_{p-1}+3y_{p-1}.y_0-3y_{0}.y_{1}+y_1.y_2\big) \in
    \bigoplus_{i=0}^{p-2} H^2(\GG,\langle h_ih_{p-1}^{-1}\rangle).
\end{multline*}
We let $\Omega_2$ be the central extension of $\H$ by $\GG$
corresponding to $y\in H^2(\GG,\H)$:
\begin{equation*}
    1\to \H \to \Omega_2 \to \GG\to 1.
\end{equation*}

The group $\Omega_2$ is a torsion-free nilpotent pro-$p$-group of
Hirsch length $2p-1$. For suitable representatives $g_i\in \Omega_2$
with respective images $\bar g_i\in \GG$, we have the relations
\begin{align*}
    [g_i,g_{i+1}] &=
    (h_{i-3}h_{p-1}^{-1})(h_{i-2}h_{p-1}^{-1})^{-3}
    (h_{i-1}h_{p-1}^{-1})^{3}
    (h_{i}h_{p-1}^{-1})^{-1}\\
    [g_i,g_j] &= 1, \qquad\qquad\qquad\qquad\qquad
    \qquad j\neq i+1\\
    [g_i,h_jh_k^{-1}] &= 1, \qquad\qquad\qquad\qquad
    \qquad\qquad \text{all\ }i,j,k\\
    [h_ih_j^{-1},h_kh_l^{-1}] &= 1, \qquad\qquad\qquad
    \qquad\qquad\qquad\text{all\ }i,j,k,l.
\end{align*}
We check that the action of $C$ on $\bar g_i$ may be extended to an
action of $C$ on the abelian pro-$p$-group $\H$ and then on
$\Omega_2$, via $\sigma(g_i)=g_{i-1}$, $\sigma(h_ih_{j}^{-1}) =
h_{i-1}h_{j-1}^{-1}$, $0\le i,j< p$.  In particular, $H^1(\H,\Fp)$
is an $\Fp C$-module isomorphic to $M_{p-1}$.

The $E_2$ page of the spectral sequence
\begin{equation*}
    H^s(\GG,H^t(\H,\Z_p)) \Rightarrow H^{s+t}(\Omega_2,\Z_p)
\end{equation*}
is generated by the anticommuting elements $y_i\in E_2^{1,0}$ and
$z_j\in E^{0,1}$, and the differential $d_2$ is given by
\begin{equation*}
    d_2^{0,1}(z_i) = -y_{i}.y_{i+1}+3y_{i+1}.y_{i+2}-3y_{i+2}.y_{i+3}
    +y_{i+3}.y_{i+4} = (\sigma-1)^3 y_i y_{i+1}.
\end{equation*}
(See \cite[Theorem~2.1.8]{NSW}.) The same holds for the $E_2$ page
of the spectral sequence with coefficients in $\Fp$
\begin{equation*}
    H^s(\GG,H^t(\H,\Fp)) \Rightarrow H^{s+t}(\Omega_2,\Fp),
\end{equation*}
and we write $\bar y_i$ and $\bar z_j$ for the reduction mod $p$ of
the classes $y_i$ and $z_j$. The picture below indicates the bases
of relevant cohomology groups on this $E_2$ page.
\[ \begin{array}{|c|c|c|c}
 & & & \\ \hline \bar z_i.\bar z_j & \bar z_i.\bar z_j.\bar y_k &  &
\\ \hline \bar z_i & \bar z_i.\bar y_j & \bar z_i.\bar y_j.\bar y_k
& \\ \hline 1 & \bar y_i & \bar y_i.\bar y_j & \bar y_i.\bar
y_j.\bar y_k \\ \hline
\end{array} \]

First, we consider $d_2^{0,1}$ on $E^{0,1}=H^1(\H,\F_p)\simeq
M_{p-1}$. The group $H^2(\GG,\Fp)\simeq \bigwedge^2 H^1(\GG,\Fp)$ is
a direct sum of free $\Z_p C$-modules $M(0,i)$, $1\le i\le (p-1)/2$,
each with $\Z_p C$-module generator $y_0.y_i$. For each $i$ in $1\le
i\le (p-1)/2$, the module $M(0,i)$ consists of the span of the
twofold products $y_j.y_{i+j}$, $0\le j<p$.  From $\im
d_2^{0,1}\subset M(0,1)$ we deduce that $\im d_2^{0,1} =
(\sigma-1)^3M(0,1)$.  Since $H^1(\H,\Fp) \simeq M_{p-1}$ we find
that $d_2^{0,1}$ has kernel $(\sigma-1)^{p-3}M_{p-1}\simeq M_{2}$,
image $M_{p-3}$, and cokernel $M_3\oplus \frac{p-3}{2} M_p$.

We may now calculate $H^1(\Omega_2,\Fp)$.  From the Five-Term Exact
Sequence we have that $H^1(\Omega_2,\Fp)$ is an extension of the
free $\Fp C$-module $H^1(\GG,\Fp)$ by $\ker d_2^{0,1} \simeq M_2$.
Since free $\Fp C$-modules are injective (see \cite[Theorem
11.2]{Ca}), we deduce that $H^1(\Omega_2,\Fp)\simeq M_2\oplus M_p$.

Next we compute the kernel of $d_2^{1,1}$. Multiplying elements in
$\ker d_2^{0,1}$ by arbitrary elements of $E_2^{1,0}$ gives elements
in $\ker d_2^{1,1}$.  Since $\ker d_2^{0,1}\simeq M_2$ and
$E_2^{1,0}\simeq M_p$, we have that $\ker d_{2}^{1,1}$ contains a
$\Fp C$-submodule isomorphic to $M_2\otimes M_p\simeq M_p\oplus
M_p$. We claim that $\ker d_2^{1,1}$ is in fact isomorphic to
$M_p\oplus M_p$, as follows. The number of cyclic $\Fp C$-summands
of an arbitrary $\Fp C$-module $M$ is equal to $\dim_{\Fp} M^C$ (see
\cite[Proposition~2]{LMS}). Moreover, $E_2^{1,1} =
H^1(\GG,\Fp)\otimes H^1(\H,\Fp)\simeq M_p\otimes M_{p-1} \simeq
(p-1) M_p$, and so $(E_2^{1,1})^C = (\sigma-1)^{p-1} E_2^{1,1}$.
Therefore it suffices to prove that the kernel $d_2^{1,1}$ on the
image of
\begin{equation*}
    N_\sigma := 1+\sigma+\dots+\sigma^{p-1} = (\sigma-1)^{p-1}
\end{equation*}
on $E_2^{1,1}$ has dimension at most $1$.

Let
\begin{equation*}
    \nu = \sum_{j=0}^{p-1}\alpha_j\sum_{k=0}^{p-1} \bar z_k.\bar
    y_{j+k} = \sum_{i,j}\alpha_{j-i}\bar z_i.\bar y_j,
\end{equation*}
be an arbitrary element of $E_2^{1,1}$, where $\alpha_i\in \Fp$.
Since $\sum_{k=0}^{p-1} \bar z_k=0$, we may add the same quantity to
each of the $\alpha_j$ and leave $\nu$ unchanged. Therefore we
assume that $\alpha_0=0$, and the remaining $\alpha_j$ are uniquely
defined. We calculate
\begin{align*}
    d_2^{1,1}(\nu) &= \sum_{i,j}\alpha_{j-i}[(\sigma-1)^3(\bar
    y_i.\bar y_{i+1})].\bar y_j \\ &= \sum_{i,j}\Big( -\alpha_{j-i}
    +3\alpha_{j-i-1} - 3\alpha_{j-i-2} +\alpha_{j-i-3}\Big) \bar
    y_i.\bar y_{i+1}.\bar y_j.
\end{align*}
If $d_2^{1,1}(\nu)=0$, the coefficient of each $\bar y_i.\bar
y_j.\bar y_k$ has to be zero. We must be careful, however, because
some monomials occur in more than one way, and monomials with
repeated indices are zero.

Recalling that $\alpha_0=0$, suppose that we are further given
$\alpha_1$ and $\alpha_2$ for an element $\nu$ in the kernel. Then
an examination of the coefficient of $\bar y_0.\bar y_1.\bar y_3$
shows that $-\alpha_0+3\alpha_1-3\alpha_2+\alpha_3=0$, and so
$\alpha_3$ is determined. Proceeding inductively on $i$, the
coefficient of $\bar y_0.\bar y_1.\bar y_i$ determines $\alpha_i$,
except for the last coefficient $\alpha_{-1}$, because $\bar
y_{-1}.\bar y_0.\bar y_1$ appears in more than one way.
Nevertheless, considering the term $\bar y_3.\bar y_4.\bar y_1$, we
find that
\begin{equation*}
    -\alpha_{-2}+3\alpha_{-1}-3\alpha_{0}+\alpha_1.
\end{equation*}
Since $\alpha_0 = 0$ we see that $\alpha_{-2}$ and $\alpha_1$
determine $\alpha_{-1}$.  This completes the proof that the kernel
of $d_2$ on $E_2^{1,1}$ is a direct sum of $2$ copies of $M_p$.
Since free $\Fp C$-modules are injective, $\ker d_2^{1,1}$ is a
direct summand of $E_2^{1,1}$ and so
$E_\infty^{1,1}=E_3^{1,1}=(p-3)M_p$.

We now turn to the computation of the kernel of $d_2^{0,2}$. We
claim that $\ker d_2^{0,2} =\bigwedge^2 \ker d_2^{0,1} \simeq M_1$.
Since $\bar z_i$, $0\le i<p-1$, form a basis of $E_2^{0,1}$, the
monomials $\bar z_i.\bar z_j$, $0\le i < j < p-1$ form a basis for
$E_2^{0,2}$. Similarly, the monomials $\bar z_k.\bar y_i.\bar y_j$
for $0\le k<p-1$ and $0 \le i < j \le p-1$ form a basis for
$E_2^{2,1}$. We have
\begin{equation*}
    d_2^{0,2}(\bar z_i.\bar z_j) = \bar z_i.[(\sigma-1)^3
    (\bar y_j.\bar y_{j+1})] - \bar z_j.[(\sigma-1)^3
    (\bar y_i.\bar y_{i+1})], \quad i<j.
\end{equation*}

Suppose that $\gamma =\sum_{i<j} \beta_{i,j}\bar z_i.\bar z_j$ lies
in the kernel. We shall show that $\beta_{12}$ determines the
remaining $\beta_{ij}$. The coefficient $\lambda_{1,3,4}$ of $\bar
z_1.\bar y_3. \bar y_4$ in $d_2(\gamma)$ is determined by the images
of monomials $\bar z_1.\bar z_j$ for $j>1$ since monomials $\bar
z_i.\bar z_j$ with $i<j$ form a basis; in particular,
$\lambda_{1,3,4}$ is determined by $d_2^{0,2}(\beta_{12} \bar
z_1.\bar z_2 + \beta_{13} \bar z_1.\bar z_3)$.  Hence $\beta_{13}$
is determined by $\beta_{12}$ and $\lambda_{1,3,4}$.

For the remaining $\beta_{ij}$, we proceed by induction.  Order the
monomials $\bar z_i.\bar z_j$, $i<j$, lexicographically, and assume
that $\beta_{kl}$ is already determined for all $\bar z_k.\bar z_l <
\bar z_i. \bar z_j$ for some fixed $0 \le i < j < p-1$. Now consider
the coefficient $\lambda_{i,j,j+1}$ of $\bar z_i.\bar y_j.\bar
y_{j+1}$ in $d_2^{0,2}(\gamma)$. The summands $\beta_{kl} \bar
z_k.\bar z_l$ of $\gamma$ contributing to $\lambda_{i,j,j+1}$
include $\beta_{ij} \bar z_i.\bar z_j$ and various $\beta_{kl} \bar
z_k.\bar z_l$ with $\bar z_k.\bar z_l < \bar z_i.\bar z_j$. By
induction, then, $\beta_{ij}$ is determined by $\lambda_{i,j,j+1}$
and $\beta_{kl}$ for $\bar z_k.\bar z_l<\bar z_i.z_j$.  We conclude
that $\ker d_2^{0,2}$ is 1-dimensional and so is $\bigwedge^2 \ker
d_2^{0,1}\simeq M_1$.

Next, we observe that $E_3^{0,2}$ consists of products of elements
of $E_3^{0,1}$, and hence $d_3^{0,2}\colon E_3^{0,2}\to E_3^{3,0}$
is the zero map. We conclude from these calculations above that
$E_\infty^{i,j}=E_3^{i,j}$ for $i+j\le 2$, and this portion of the
$E_\infty$ page is as follows:
\[ \renewcommand{\arraystretch}{1.2}
\begin{array}{|c|c|c|c}
& & & \\ \hline M_1 & & & \\ \hline M_2 & 2M_p & & \\
\hline M_1 & M_p & M_3 \oplus \frac{p-3}{2}M_p & \\ \hline
\end{array} \]

We are now ready to assert that $H^2(\Omega_2,\Fp)$ is the direct
sum of the $E_\infty^{i,j}$ with $i+j=2$. To see this, we use an
automorphism $\phi$ of $\Omega_2$ such that $\phi(g_i)=g_i^2$ for
$0\le i\le p-1$ and $\phi(h_ih_{p-1}^{-1})=(h_ih_{p-1}^{-1})^4$ for
$0\le i < p-1$. That $\phi$ exists one can see from the fact that
the commutator map in a nilpotent group of class $2$ is a bilinear
map. We see that $\phi$ commutes with the action of $\sigma$ and the
induced action on the spectral sequence satisfies $\phi^*(\bar y_i)
= 2^{-1}\bar y_i$ and $\phi^*(\bar z_i)=4^{-1}\bar z_i$, where
$a^{-1}$ denotes the inverse of $a$ modulo $p$.   We calculate that
$\phi$ acts on $E_\infty^{0,2}$ by multiplication by $16^{-1}$, on
$E_\infty^{1,1}$ by multiplication by $8^{-1}$, and on
$E_\infty^{2,0}$ by multiplication by $4^{-1}$. Since $16^{-1}$,
$8^{-1}$ and $4^{-1}$ are distinct modulo $p$, we obtain that
$H^2(\Omega_2,\Fp)$ decomposes as a direct sum of the eigenspaces of
$\phi$ and so $H^2(\Omega_2,\Fp)\simeq M_1 \oplus M_3 \oplus
\frac{p+1}{2} M_p$.

Finally, since $E_\infty^{0,2}=\bigwedge^2 E_{0,1}$,
$E_\infty^{1,1}=E_{0,1}.E_{1,0}$, and $E_\infty^{2,0}$ is a quotient
of $E_2^{2,0}=\bigwedge^2 E_{0,1}$, we deduce that
$H^2(\Omega_2,\Fp)$ is decomposable.

\begin{remark*}
    The relations among the $g_i$ and $h_ih_{p-1}^{-1}$ given in the
    construction of $\Omega_2$ are minimal since the differential
    $d_2^{0,1}$ on the $E_2$ page of the spectral sequence with
    coefficients in $\Z_p$, $H^s(\GG,H^t(\H,\Z_p)) \Rightarrow
    H^{s+t}(\Omega_2,\Z_p)$, is injective, as follows.  The
    differential $d_2^{0,1}$ is given by $d_2^{0,1}(z_i) =
    (\sigma-1)^3(y_i.y_{i+1})$. Observe that $H^1(\H,\Z_p) \simeq
    E_2^{0,1}$ is isomorphic as a $\Z_p C$-module with $\Z_p[\xi_p]$,
    the ring of integers in $\Q_p[\xi_p]$, and that $\sigma$ acts on
    $\Z_p[\xi_p]$ as multiplication by $\xi_p$. We deduce that
    multiplication by $(\sigma-1)^3$ is injective. Moreover, the map
    \begin{equation*}
        e \colon E_2^{0,1} \to E_2^{2,0}, \quad\quad e(z_i) = y_i
        y_{i+1}, \quad\quad 0\le i< p-1,
    \end{equation*}
    is injective. Hence $d_2^{0,1}$, as the composition of injective maps,
    is injective.
\end{remark*}

\section*{Acknowledgement}

We thank the referee of \cite{LMS} for suggesting that we illustrate
\cite[Theorem 1]{LMS} by constructing examples of pro-$p$-groups
which are not absolute Galois groups, for this appendix is the
result.  The third author would also like to thank Sunil Chebolu,
Ajneet Dhillon and Vahid Shirbisheh for stimulating discussions
related to this paper.

\end{document}